\documentclass[10pt]{amsart}

\usepackage[latin2]{inputenc}
\usepackage{amsmath}
\usepackage{graphicx}
\usepackage{amssymb}
\usepackage{bm}
\usepackage{bm}
\usepackage{esint}
\usepackage{color}
\usepackage{amsthm}
\usepackage{epsfig}
\usepackage{enumitem}
\usepackage{mathtools}
\usepackage{esvect}
\usepackage{tikz}
\usepackage{mathrsfs}
\usetikzlibrary{calc,intersections,through,backgrounds,patterns,arrows.meta, math,through}
\usepackage[english]{babel}
\usepackage[linktocpage=true,colorlinks=true,linkcolor=Blue,citecolor=Green]{hyperref}
\usepackage{pgfplots}
\usepackage{pgfmath}
\pgfplotsset{compat=1.17}
\usetikzlibrary{arrows.meta}

\newtheorem{theorem}{Theorem}
\newtheorem{proposition}[theorem]{Proposition}
\newtheorem{lemma}[theorem]{Lemma}

\newtheorem*{theorem*}{Theorem}

\theoremstyle{definition}
\newtheorem{definition}{Definition}
\newtheorem{remark}[definition]{Remark}

\def\XXint#1#2#3{{\setbox0=\hbox{$#1{#2#3}{\int}$ }
\vcenter{\hbox{$#2#3$ }}\kern-.6\wd0}}

\definecolor{Yellow}{rgb}{0.95,0.9,0.0} 
\definecolor{Red}{rgb}{0.8,0.1,0.1}
\definecolor{Green}{rgb}{0.1,0.65,0.2}
\definecolor{Blue}{rgb}{0.1,0.1,0.8}
\definecolor{Purple}{rgb}{0.7,0.1,0.7}
\definecolor{Grey}{rgb}{0.6,0.6,0.6}

\definecolor{YELLOW}{rgb}{0.95,0.9,0.0} 
\definecolor{RED}{rgb}{0.8,0.1,0.1}
\definecolor{GREEN}{rgb}{0.25,0.65,0.1}
\definecolor{BLUE}{rgb}{0.1,0.1,0.8}
\definecolor{PURPLE}{rgb}{0.7,0.1,0.7}

\newcommand{\supp}{\operatorname{supp}}

\DeclareMathOperator*{\essinf}{ess\,inf}

\newcommand{\Rd}[1][d]{{\mathbb{R}^{#1}}}
\allowdisplaybreaks[4]

\newcommand{\R}{\mathbb{R}}
\newcommand{\dH}{\mathrm{d}\mathcal{H}^{d-1}}

\newcommand{\dx}{\mathrm{d}x}
\newcommand{\dt}{\mathrm{d}t}

\begin{document}

\title[The Verigin problem with phase transition] 
{The Verigin problem with phase transition as a Wasserstein flow}

\author{Anna Kubin}
\address{Institut f\"ur Analysis und Scientific Computing, Technische Universit\"at Wien, Wiedner Hauptstrasse 8-10, 1040 Vienna, Austria}
\email{anna.kubin@tuwien.ac.at}

\author{Tim Laux}
\address{Institut f\"{u}r Matheamtik \& Interdisziplin\"{a}res Zentrum f\"{u}r Wissenschaftliches Rechnen, Universit{\"a}t Heidelberg,  Im Neuenheimer Feld 205, 69120 Heidelberg, Germany.}
\email{tim.laux@math.uni-heidelberg.de}

\author{Alice Marveggio}
\address{Institute for Applied Mathematics, Universit{\"a}t Bonn, Endenicher Allee 62, 53115 Bonn, Germany}
\email{alice.marveggio@hcm.uni-bonn.de}

\begin{abstract}
	We study the modeling of a compressible two-phase flow in a porous medium. The governing free boundary problem is known as the Verigin problem with phase transition.
	We introduce a novel variational framework to construct weak solutions.
    Our approach reveals the gradient-flow structure of the system by adopting a minimizing movement scheme using the Wasserstein distance.
	We prove the convergence of the scheme, obtaining ``relaxed" distributional solutions in the limit that satisfy an optimal energy-dissipation rate. 
	Under the additional assumptions that $d \geq 3$ and that the discrete mass densities are uniformly Muckenhoupt weights, we show that the limit is the characteristic function of a set of finite perimeter in the region where there is no vacuum.

    \vspace{5pt}

    \vspace{5pt}
    \noindent
    \textbf{Mathematical Subject Classification (MSC) 2020:} 
    35R35, 
    49Q22, 
    35R37, 
    76D27, 
    76S05, 
    35A01, 
    35D30. 
\end{abstract}


\maketitle

\section{Introduction}
The Verigin problem models a compressible two-phase flow in a porous medium.
This situation occurs, for instance, when a porous medium containing a compressible fluid is flooded with another compressible liquid to displace the fluid under isothermal conditions.
Additionally, phase transition may occur and be considered into the model.
The resulting system corresponds to the Verigin problem with phase transition, a free boundary problem that is the compressible analogue of the Muskat problem (which considers incompressible phases).

 For a smooth domain $\Omega \subset \R^d$ and a final time horizon $T>0$, denoting by $\rho: \Omega \times [0,T) \rightarrow [0,\infty)$ the total density of the two-phase fluid, by $\chi: \Omega \times [0,T) \rightarrow  \{0,1\}$ the indicator function of the set occupied by one of the two fluids, and by $\Gamma :=  \partial \{ \chi =1\}$ the phase interface, the system of equations reads
 \begin{align}\label{eq:smooth_system_intro}
 	\begin{cases}
 		\partial_t \rho + \operatorname{div}(\rho \mathbf{u})=  0  & \; \text{ in } \Omega , \\
 		[\![\rho]\!]V_\Gamma =	[\![ \rho  \mathbf{u} \cdot \nu_\Gamma ]\!] & \;  \text{ along } \Gamma , 
 	\end{cases}
 	\quad \quad 
 	\begin{cases}  
 		\rho \mathbf{u}= - \nabla \pi    & \; \text{ in } \Omega ,  \\
 		[\![\pi]\!]= \sigma H_{\Gamma}  &\;  \text{ along } \Gamma.
 	\end{cases}
 \end{align}    
 Here, $\mathbf{u}:\Omega \times [0,T) \rightarrow \R^d$ is the velocity vector field, $\pi = \pi(\rho, \chi)$ is the pressure, and the constant $\sigma >0$ is the surface tension. Moreover,
 $	[\![\cdot ]\!] $ denotes the jump across the interface $\Gamma$, $V_\Gamma$ is the velocity of $\Gamma$ in the direction of the normal $\nu_\Gamma$, and $H_\Gamma$ is the (scalar) mean curvature of $\Gamma$.
 The system \eqref{eq:smooth_system_intro} is endowed with the natural boundary conditions
 \begin{align} \label{eq:boundcond_intro}
 	\nabla\rho \cdot \nu_{\Omega} = 0  \; \; \;   \text{ along } \partial \Omega, \quad 
 	\nu_\Gamma \cdot \nu_{\Omega} = 0 \; \; \; \text{ along } \partial\Omega \cap  \Gamma,
 \end{align}
 and with initial data $\chi(\cdot , 0)= \chi_0 $ and $\rho(\cdot,0)=\rho_0$.
 To the best of the authors' knowledge, the only mathematical results available in the literature on this model concern the theory of classical smooth solutions (see \cite{PrussSim,PrussSimWielke} and the references therein). In particular, short-time existence, stability, as well as convergence to equilibrium of solutions that do not develop singularities are studied in \cite{PrussSim}.

In this manuscript we establish the first long-time existence result for weak solutions to the Verigin problem with phase transition. We construct distributional solutions via an implicit time discretization scheme exploiting the gradient flow structure of the system \eqref{eq:smooth_system_intro}. 
The underlying energy of the system consists of three terms, namely: the total surface area of the interface between the fluids, the mass-specific free energy, and a coupling term favoring the phase-transition. More precisely, the total energy is (up to an additive constant) given  by
\begin{align} \label{eq:energy_intro}
	\mathcal{E}(\rho, \chi) = \sigma \int_{\Gamma} 1 \, \mathrm{d}\mathcal{H}^{d-1}  +   \int_{\Omega}    f(\rho)+  \chi \big( \lambda   - \rho\big )  \, \mathrm{d}x \,.
\end{align}
where $ \lambda>0$ is a constant and the mass-specific free energy density $f$ is given by
\begin{align*} 
	f(\rho) = \begin{cases}
		\rho \log(\rho) \,, \quad  \text{if } m=1,
		\\
		\frac{1}{m-1}	 \rho^m	 \,, \quad  \text{if }  m>1.
	\end{cases}
\end{align*}
We note that the dissipation of the energy is given by 
\[
\frac{\mathrm{d}}{\mathrm{d}t } \mathcal{E} (\rho, \chi)=  	- \int_{\Omega} \rho |\mathbf{u}|^2 \, \mathrm{d} x.
\] 
Our first insight is that, just like in the case of diffusion equations (cf.\ the seminal work of Jordan, Kinderlehrer, and Otto~\cite{JKO}), this energy dissipation relation indicates that the system can be viewed as a gradient flow with respect to the Wasserstein distance $W_2$, which plays a crucial role in optimal transport. 
Indeed, given probability densities $\rho_0$ and $\rho_T$ on $\Omega$, the Benamou-Brenier formula reads as
$$ W^2_2(\rho_0, \rho_T) = \inf_{(\rho, \mathbf{u})} \int_0^T \int_\Omega \rho(x,t) |\mathbf{u}(x, t)|^2\, \dx 	\dt$$
subject to the continuity equation $	\partial_t \rho + \operatorname{div}(\rho \mathbf{u})=  0 $ and the boundary conditions $\rho(\cdot, 0)= \rho_0$ and $\rho(\cdot, T)= \rho_T$.
This clearly suggests that ``the mass density $\rho$ evolves according to the gradient flux of the energy \eqref{eq:energy_intro} with respect to the 2-Wasserstein metric".
This idea has found many applications in the context of partial differential equations; however, the literature on Wasserstein gradient flows in the context of free boundary problems is much smaller and still evolving, see e.g.~\cite{Otto2, Chambolle-Laux, JKM_Muskat}.

A key feature of our problem lies in 
the interplay between the diffusion dynamics of the mass density $\rho$ and the quasi-static evolution of the free boundary $\Gamma$. The system \eqref{eq:smooth_system_intro} describes a two-phase fluid flow of mass density $\rho$. The evolution of the interface $\Gamma$  is governed by the normal flux of the fluid through $\Gamma$, including mass transfer due to phase change. 
 Since $\mathbf{u}$ is determined by Darcy's law in the bulk, we can view it together with the transport equation for $\rho$ as a single diffusion equation. The jump of the pressure $\pi$ across $\Gamma$ is related to the mean curvature by the Young-Laplace law, which provides the dynamic boundary condition linking the interface geometry to the bulk flow.
Similar quasi-static dynamics underlie the Stefan problem with Gibbs-Thomson law studied in \cite{Luckhaus1}, where the interface velocity is determined by the heat flux across the interface, while the jump of the temperature across the interface is related to the mean curvature according to the Gibbs-Thomson relation. In the spirit of Luckhaus' work \cite{Luckhaus1}, we show that our discrete flow $(\rho^h, \chi^h)$ satisfies an almost minimizing property, which allows us to deduce the $L_t^1$-convergence of the perimeters as the time discretization parameter $h>0$ vanishes. This property is crucial for our analysis: it leads to an existence result for ``relaxed" distributional (or $BV$) solutions to the Verigin problem with phase transition without assuming a priori perimeter convergence. Such an assumption was originally introduced in the seminal work \cite{LuckStur} and has subsequently largely appeared in the literature under the name of ``energy convergence hypothesis". We emphasize that this hypothesis is not required in the framework of varifold solutions (see, e.g., \cite{ATW,Roger_MS,Roger_SP,HenselStinson}), where one does not need to ensure unit multiplicity of the limiting interface.

To strengthen our existence result and lift the word ``relaxed", we need to show that the limit as $h \rightarrow 0$ of $\chi^h$ is indeed the indicator function of a set of finite perimeter. To this aim, we improve the convergence properties of $\chi^h$ by establishing strong convergence in time. This can be achieved by transferring the time compactness from the mass density $\rho^h$ to $\chi^h$ via a crucial stability estimate. However, this estimate can be proved only assuming that the functions $(\rho^h \wedge 1)^\frac{d}{d-1}$ are uniformly Muckenhoupt weights. Under this condition, and the technical assumption $d \geq 3$, the limit $\chi$ is the indicator function of a set of finite perimeter in the region where $\rho$ is strictly positive. This limitation is natural: when 
$\rho$ vanishes, the notion of distinct phases loses its meaning, as one is effectively in vacuum.

Finally, we briefly situate our result within the current literature. Only a few works interpret free boundary problems as Wasserstein gradient flows.
A Wasserstein gradient-flow structure of the Hele-Shaw flow has been observed in \cite{Otto2} and later exploited in \cite{Chambolle-Laux} to establish the existence of distributional solutions via a minimizing movement scheme. Distributional solutions to the incompressible analogue of our model, i.e., the Muskat problem, were constructed in \cite{JKM_Muskat}, conditional on an energy convergence assumption. The Verigin problem with phase transition presents additional challenges beyond these settings, notably the joint minimization in two variables and the possible vanishing of the mass density. On the other hand, we show that for the Verigin problem with phase transition, the convergence of the approximation scheme is better than in the case of many other free boundary problems. In particular we can \emph{prove} the perimeter convergence, which oftentimes is an \emph{assumption}~\cite{LuckStur} or holds only under specific assumptions on the geometry~\cite{DePhilippisLaux}.

The rest of the paper is organized as follows. We derive the model in Section \ref{sec:model}. We present the precise mathematical setting, the minimization problem, and the main results in Section \ref{sec:results}. Section \ref{sec:discreteprob} concerns the properties of the discrete problem. Compactness properties of the discrete flow are shown in Section \ref{sec:compactness}.
In Section \ref{sec:existence}, we prove the existence of ``relaxed" distributional solutions to the Verigin problem with phase transition. Finally, in Section \ref{sec:strongconv}, we improve the convergence of the indicator functions of the discrete sets under the Muckenhoupt assumption.

\vspace{2mm}

\paragraph{\bfseries Notation.}
We denote by $d$ an integer number larger than or equal to 2.
Let $\Omega \subset \mathbb{R}^d$ be an open set, possibly $\Omega = \R^d$, with boundary $\partial \Omega$ and outer normal vector field $\nu_\Omega$.
	For any measurable set $E \subset \Omega$ we denote by $\chi_E$ its characteristic function and by $\nu_{E}$ its outer normal vector field. 
Given a set $E \subset \Omega$ and $g:\Omega \to \R$ a measurable function, we denote by $[\![g]\!](x):=g_+(x)-g_-(x)$ for $x \in \partial E$, where $g_+=g|_{\Omega \setminus E}$ and $g_-= g|_E$.
We say $E$ is a set of finite perimeter in $\Omega$ if $\chi_E\in BV(\Omega;\{0,1\}),$ and we denote its distributional derivative in $\Omega$ by $\nabla \chi_E$ and its total variation measure by $|\nabla \chi_E|$. 
We denote by $P(E; \Omega):= \int_\Omega 1 \, \mathrm{d} |\nabla \chi_E|$, omitting $\Omega$ when $\Omega=\R^d$, and we use $\nu_E:=-\frac{\mathrm{d}\nabla \chi_E}{\mathrm{d} |\nabla \chi_E|}$ for the measure-theoretic outer normal of $E$. 
Using $\partial^\ast E$ as the reduced boundary of the set $E$, we have $\int_\Omega \cdot \, \mathrm{d} |\nabla \chi_E| =  \int_{\partial^\ast E \cap \Omega} \cdot \, \mathrm{d} \mathcal{H}^{d-1}$, where $\mathcal{H}^{d-1}$ is the Hausdorff measure.
Finally, for $p >1$, we denote the Muckenhoupt class by $A_p$ (see Appendix \ref{sec:Mucken}).

\section{Derivation of the model}\label{sec:model}
We aim to model a compressible two-phase flow in a porous medium, allowing for phase transition between the two phases. 
Our derivation of the model is inspired by the works of Pr\"uss, Simonett, and Wielke \cite{PrussSim,PrussSimWielke}.
The resulting model is the Verigin problem with phase transition, which is the compressible analogue of the Muskat problem.

Consider a smooth bounded domain $\Omega \subset \Rd$ with boundary $\partial \Omega$ and outer normal vector field $\nu_{\Omega}$.
Let $\rho_+,\rho_-: \Omega \times [0, \infty) \rightarrow [0, \infty)$ be the densities of two distinct compressible fluids. We denote by  $\Omega_+ $ and by $\Omega_- $ the subdomains of $\Omega$ occupied by the two fluids, respectively.
Let $\chi: \Omega \times [0, \infty) \rightarrow \{0,1\}$ be the phase indicator function associated to $\Omega_-$, let $\Gamma:= \partial\{\chi=1\}$ be the boundary of $\Omega_-$ such that $\Omega = \Omega_+ \cup \Omega_- \cup \Gamma$, and let $\nu_\Gamma$ be the outer normal vector field of $\Gamma$. 
We can define a single density function as $\rho = \rho_- \chi + \rho_+ (1-\chi)$. 
Furthermore, we introduce a pressure field $\pi $ and a sufficiently general equation of state including all fluids and conditions of flow of practical interest (cf. \cite{Muskat}): $$\rho= \rho_0 \pi^\alpha \exp(\beta \pi )  ,$$ where $\rho_0, \alpha, \beta>0$  are constants. 
In particular, we consider  \begin{align}
	\rho_\pm= 	\begin{cases}
		\rho^\pm_0 \exp(\beta^\pm \pi ) \simeq \rho^\pm_0 (1+ \beta^\pm \pi   ) \; &\text{for a compressible liquid},\\
		\rho^\pm_{0}  \pi  \; &\text{for a gas under isothermal expansion}, \\
		\rho^\pm_{0}  \pi^\alpha, \, 0 < \alpha < 1, \; &\text{for a gas under adiabatic expansion}.
	\end{cases}
\end{align}

In order to derive the model, we recall that conservation of mass reads as
\begin{align}
	\partial_t \rho + \operatorname{div}(\rho \mathbf{u}) = 0 \quad \text{ in } \Omega \setminus \Gamma,\label{eq:conservation}  \\
	[\![\rho (\mathbf{u} \cdot \nu_\Gamma - V_\Gamma)]\!] = 0  \quad \text{ along } \Gamma , \label{eq:jumpflux} 
\end{align}
where $	[\![\cdot ]\!] $ denotes the jump across $\Gamma$, $V_\Gamma$ is the normal velocity of the interface $\Gamma$ in the direction of $\nu_\Gamma$.
Additionally, we require $\mathbf{u} \cdot \nu_{\Omega} = 0$, so that the flow induced by the velocity field $\mathbf{u}$ does not leave $\Omega$.
The jump condition in~\eqref{eq:jumpflux} shows that the phase flux $j_\Gamma := \rho (\mathbf{u} \cdot \nu_\Gamma - V_\Gamma)$ is uniquely defined on $\Gamma$. Moreover, the jump condition~\eqref{eq:jumpflux} can be rewritten as
\begin{align} 
[\![\rho]\!]V_\Gamma = 	[\![\rho \mathbf{u} \cdot \nu_\Gamma ]\!] , \quad  \frac{1}{[\![\rho]\!]}j_\Gamma = [\![ \mathbf{u} \cdot \nu_\Gamma ]\!]. 
\end{align}
On the interface we impose the Laplace-Young law and the ninety-degree angle condition
\begin{align} \label{eq:jumppi}
	[\![\pi]\!] = \sigma H_\Gamma & \quad \text{ along } \Gamma , \\
	\nu_\Gamma \cdot 	\nu_{ \Omega} = 0 &\quad \text{ along } \partial\Omega \cap  \Gamma, 
\end{align}
where the constant $\sigma >0$ denotes the surface tension and $H_\Gamma$ the mean curvature of $\Gamma$.

The total energy of the problem is given by
\begin{align} \label{eq:modelenergy}
	\mathcal{E}(\rho, \chi) = \sigma \mathcal{H}^{d-1}(\Gamma)+   \int_{\Omega} \rho \psi(\rho)  \, \mathrm{d}x, 
\end{align}
where $\psi$ is the mass-specific free energy density, which also depends on the phases (i.e., liquid or gas). 
By using Maxwell's law
\begin{align} \label{eq:modelpi}
	\pi(\rho)= \rho^2 \psi^\prime(\rho),
\end{align}
and the equations of state for $\pi$ for both phases, we deduce 
\begin{align}
	\psi^\prime(\rho )  = \begin{cases}
		\frac{\rho/\rho^\pm_0 -1}{\beta_\pm\rho^2}   \quad &\text{for a liquid}
		\\
			\frac{1}{\rho_0^\pm \rho} \quad  &\text{for an  isothermal gas} \\
			\frac{\rho^{1/\alpha -2}}{(\rho_0^\pm)^{1/\alpha}} \quad  &\text{for an adiabatic gas}
	\end{cases}\quad \text{ in } \Omega_\pm.
\end{align}

By first computing $\psi$, one can rewrite the total energy~\eqref{eq:modelenergy} as
\begin{align} \label{eq:modelgammaenergy}
	&\mathcal{E}(\rho, \chi) \\ \notag 
	&= \sigma \int_{\Gamma} 1 \, \mathrm{d}\mathcal{H}^{d-1}  +   \int_{\Omega}   (\gamma_1 \chi +   \gamma_2 )  f(\rho) +  (\gamma_3 \chi  +  \gamma_ 4 ) + (\gamma_5 \chi + \gamma_6 )\rho \, \mathrm{d}x,
\end{align}
where $\{\gamma_i\}_{i \in \{1,...,6\}}$ are constants and $f(\rho)$ is either of the form $\rho \log(\rho)$ or $\rho^{1/\alpha}$.
In particular, we have
\begin{align}
\psi(\rho, \chi  ) & = 
(\gamma_1 \chi +   \gamma_2 )   \frac{f(\rho)}{\rho} +  (\gamma_3 \chi  +  \gamma_ 4 ) \frac{1}{\rho}+ (\gamma_5 \chi + \gamma_6 )	, \\
			\pi(\rho , \chi) &= 
		(  \gamma_1 \chi +   \gamma_2 ) (\rho f'(\rho) - f(\rho))- (\gamma_3 \chi  +  \gamma_ 4 ). 
\end{align}

To obtain a closed model allowing for phase transition (corresponding to $j_\Gamma \neq 0$), we impose Darcy's law for the velocity
\begin{align} \label{eq:modeldarcy}
	\mathbf{u} = - \frac{k}{\mu} \nabla \pi , 
\end{align}
where $k>0$ is the permeability of the porous medium and $\mu=\mu(\rho)>0$ denotes the viscosity of the fluids, 
and
\begin{align} \label{eq:jumpvarphi}
[\![\psi(\rho, \chi ) + \pi(\rho, \chi )/\rho ]\!]=0 \quad \text{ along } \Gamma.
\end{align}
In particular, we assume
$\mu(\rho) ={\rho}$ in~\eqref{eq:modeldarcy}, and note that~\eqref{eq:jumpvarphi} yields
\begin{align*}
	[\![(\gamma_1 \chi +   \gamma_2 )   \log(\rho) + (\gamma_1 + \gamma_5) \chi ]\!] =0  \quad \text{ along } \Gamma.
\end{align*}
Moreover, we observe that  the chemical potential of the fluid is given by
\begin{align} \label{eq:modelphi}
\phi(\rho, \chi ) &= \partial_\rho (\rho \psi(\rho, \chi ) ) = \psi(\rho, \chi ) +  \pi(\rho, \chi )/\rho 
\\ 
&=  \notag
(\gamma_1 \chi +   \gamma_2 )  f'(\rho) +  \gamma_5\chi + \gamma_6 
, \end{align}
hence~\eqref{eq:jumpvarphi} corresponds to 
\begin{align} \label{eq:jumpphi}
    [\![\phi(\rho, \chi )]\!]=0 \quad \text{ along } \Gamma.
\end{align}

As a consequence of the modeling above, our problem reads as the following system of equations
\begin{subequations}
\begin{align} 
\label{eq:model1}
	\partial_t \rho &= k \Delta \rho && \quad \text{ in } \Omega \setminus \Gamma,  \\ \label{eq:model2}
	[\![\pi(\rho, \chi )]\!]&= \sigma H  &&\quad \text{ along } \Gamma ,\\ \label{eq:model3}
		[\![\phi(\rho, \chi ) ]\!]&=0 &&\quad \text{ along } \Gamma , \\
		[\![\rho]\!]V_\Gamma +	[\![\nabla \rho  \cdot \nu_\Gamma ]\!] &=0 &&\quad \text{ along } \Gamma , \label{eq:model5}
\end{align}
\end{subequations}
endowed with the natural boundary conditions
\begin{align*}  
       \nabla\rho \cdot \nu_{\Omega} &= 0  \; \; \;   \text{ along } \partial \Omega, \\  
        \nu_\Gamma \cdot \nu_{\Omega} &= 0 \; \; \; \text{ along } \partial\Omega \cap  \Gamma,
\end{align*}
and with initial data $\Gamma(0)= \Gamma_0 $ and $\rho(0)=\rho_0$.
Observe that, choosing $\pi$ as independent variable and $\rho$ constant in the phases, the system above reduces to another geometric equation, namely the Muskat flow with phase transition.

Furthermore, one can prove the time-dissipation of the total energy~\eqref{eq:modelenergy}:
\begin{align} \label{eq:modelendissip}
	\frac{\mathrm{d}}{\mathrm{d}t } \mathcal{E} (\rho, \chi)=   	- \int_{\Omega} \frac{k}{\rho} |\nabla \pi|^2 \,   \mathrm{d} x = 	- \int_{\Omega} \frac{\rho}{k} |\mathbf{u}|^2 \, \mathrm{d} x,
\end{align}
which follows from
\begin{align*}
	\frac{\mathrm{d}}{\mathrm{d}t } \bigg(  \sigma \int_{\Gamma} 1 \, \mathrm{d}\mathcal{H}^{d-1} \bigg) &= - \sigma \int_{\Gamma} H_\Gamma V_\Gamma \, \mathrm{d}\mathcal{H}^{d-1} ,
    \\
	\frac{\mathrm{d}}{\mathrm{d}t } \int_{\Omega} \rho \psi (\rho, \chi)  \, \dx 
	 &= \int_{\Omega \setminus \Gamma} \partial_\rho (\rho \psi(\rho)) \partial_t \rho  \, \dx - \int_{\Gamma} [\![\rho \psi(\rho,\chi)]\!]V_\Gamma \, \dH .
\end{align*}
and the relations above. Indeed, we  can compute
\begin{align*}
    \frac{\mathrm{d}}{\mathrm{d}t } \mathcal{E} (\rho, \chi) &= \int_{\Omega \setminus \Gamma} \partial_\rho (\rho \psi(\rho)) \partial_t \rho  \, \dx - \int_{\Gamma} ([\![\rho \psi(\rho,\chi)]\!]+ \sigma H_\Gamma )V_\Gamma \, \dH\\
    &= \int_{\Omega \setminus \Gamma} \rho  \mathbf{u} \cdot \nabla \phi    \, \dx 
    + \int_{\Gamma} [\![ \rho \mathbf{u} \cdot \nu_\Gamma   ]\!] \phi \, \dH
     - \int_{\Gamma} ([\![\rho \psi(\rho,\chi)]\!]+ \sigma H_\Gamma )V_\Gamma \, \dH\\
     &= \int_{\Omega \setminus \Gamma} \rho  \mathbf{u} \cdot \nabla \phi    \, \dx 
    + \int_{\Gamma} [\![ \rho    ]\!] \phi V_\Gamma \, \dH
     - \int_{\Gamma} ([\![\rho \psi + \pi (\rho, \chi) ]\!] V_\Gamma \, \dH
\end{align*}
where first we inserted~\eqref{eq:conservation}, integrated by parts, and recalled~\eqref{eq:jumpphi}, then we used~\eqref{eq:jumpflux} and~\eqref{eq:jumppi}. Notice that the sum of the last two integrals in the last line is zero due to~\eqref{eq:modelphi} and~\eqref{eq:jumpphi}. Finally,  from~\eqref{eq:modelphi} we can deduce
\begin{align*}
\nabla \phi &= \frac{1}{\rho}\nabla (\rho \psi - \pi) -\frac{1}{\rho^2}  (\rho \psi - \pi) \nabla \rho \\
&=  - \nabla  \psi - \frac{1}{\rho} \nabla \pi + \frac{1}{\rho^2} \pi \nabla \rho = - \frac{1}{\rho} \nabla \pi,
\end{align*}
where we used~\eqref{eq:modelpi}. The energy dissipation relation~\eqref{eq:modelendissip} then follows from~\eqref{eq:modeldarcy}.

\subsection{Our setting} 
By choosing $\gamma_1 = \gamma_4  =0$, $\gamma_2=1$, $\gamma_5 = -1$, $\gamma_3 = \lambda$,  and $\gamma_6 = c_{\Omega}>0$ in~\eqref{eq:modelgammaenergy}, the total energy reads as
\begin{align*} 
	\mathcal{E}(\rho, \chi) = \sigma \int_{\Gamma} 1 \, \mathrm{d}\mathcal{H}^{d-1}  +   \int_{\Omega}    f(\rho)+  \chi \Big( \lambda   - \rho\Big )  \, \mathrm{d}x \,  + \, c_{\Omega}\,,
\end{align*}
where
\begin{align*} 
	f(\rho) = \begin{cases}
		 \rho \log(\rho) \,, \quad  \text{if } m=1,
		 \\
	\frac{1}{m-1}	 \rho^m	 \,, \quad  \text{if }  m>1.
	\end{cases}
\end{align*}
In particular, we consider
\begin{align*}
		\psi(\rho, \chi  )= 
	   \frac{f(\rho)}{\rho} +  \lambda  \chi   \frac{1}{\rho}-  \chi  	, \quad 
	\pi(\rho , \chi) = 
	\rho^m -   \lambda  \chi   , \quad 
	\phi (\rho, \chi )  =      f'(\rho) - \chi +    1. 
\end{align*}

Under these assumptions, the Verigin problem with phase transition~\eqref{eq:model1}-\eqref{eq:model5} can be rewritten as
\begin{align}\label{eq:smooth_system}
\begin{cases}
	\partial_t \rho + \operatorname{div}(\rho \mathbf{u})=  0  & \; \text{ in } \Omega , \\
    [\![\rho]\!]V_\Gamma =	[\![ \rho  \mathbf{u} \cdot \nu_\Gamma ]\!] & \;  \text{ along } \Gamma , 
\end{cases}
\quad \quad 
\begin{cases}  
    \rho \mathbf{u}= - \nabla \pi    & \; \text{ in } \Omega ,  \\
	[\![\pi]\!]= \sigma H  &\;  \text{ along } \Gamma,
\end{cases}
\end{align}    
endowed with boundary conditions
\begin{align} \label{eq:boundcond1}
       \nabla\rho \cdot \nu_{\Omega} &= 0  \; \; \;   \text{ along } \partial \Omega, \\ \label{eq:boundcond2} 
        \nu_\Gamma \cdot \nu_{\Omega} &= 0 \; \; \; \text{ along } \partial\Omega \cap  \Gamma,
\end{align}
and with initial data $\chi(\cdot , 0)= \chi_0 $ and $\rho(\cdot,0)=\rho_0$.
Moreover, the total energy dissipation~\eqref{eq:modelenergy} reads as 
\begin{align*} 
	\frac{\mathrm{d}}{\mathrm{d}t } \mathcal{E} (\rho, \chi)=   	- \int_{\Omega} \frac{1}{\rho} |\nabla \rho|^2 \,   \mathrm{d} x = 	- \int_{\Omega} \rho |\mathbf{u}|^2 \, \mathrm{d} x.
\end{align*}

\section{Minimization problem and main results}\label{sec:results}
Let $\Omega \subset \R^d$ be an open set and $\lambda \in (0,\infty)$.
For any measurable set $E \subset \Omega$ and any measurable function $\rho:\Omega \to [0,+\infty)$, we consider
\begin{align}\label{eq:ourenergy}
\mathcal{E}(\rho, \chi_E) = \int_{\Omega} 1 \, \mathrm{d}|\nabla \chi_E |  +   \int_{\Omega}    f(\rho) \,\dx+  \int_{\Omega}\chi_E \Big( \lambda   - \rho\Big )  \, \mathrm{d}x + c_{\Omega},
\end{align}
where $c_{\Omega}>0$ and
\begin{align*} 
	f(\rho) = \begin{cases}
		 \rho \log(\rho) \,, \quad  \text{if } m=1,
		 \\
	\frac{1}{m-1}	 \rho^m	 \,, \quad  \text{if }  m>1.
	\end{cases}
\end{align*}
Moreover, we define the set
\begin{equation}
    K:= \Big \{\rho:\Omega \to (0,\infty): \rho \text{ measurable, } \int_\Omega \rho \, \dx = 1, \, \,  M(\rho)<\infty \Big \}
\end{equation}
where $M(\rho):=\int_{\Omega} |x|^2\rho(x) \,\dx$.

\begin{remark}
For convenience, we choose the constant $c_{\Omega}>0$
such that the infimum of \eqref{eq:ourenergy} over $\rho \in K$ is greater than or equal to zero for every set $E$.
Of course, the constant $c_{\Omega}$ does not affect the dynamics.
\end{remark}
\begin{remark}
We exclude the case $\lambda = 0$, since the minimization of \eqref{eq:ourenergy} with respect to $\chi_E$ becomes trivial. Indeed, if $\lambda= 0$, then then whole domain $\Omega$ is the minimizer.
For $\Omega=\R^d$, the minimization is not trivial for any $\lambda>0$, as $\mathcal{E}(\rho,\chi_\Omega)=+\infty$. 
On the other hand, when $|\Omega|<\infty$, choosing for instance $\lambda=\frac{1}{|\Omega|}$, the term $\lambda-\rho$ changes sign in $\Omega$, preventing the minimization from being trivial.
\end{remark}

Fix a time discretization step $h>0$. Given $E_0 \subset \Omega$  and $\rho_0 \in L^1(\Omega)$ such that $\int_\Omega \rho_0 \,\dx=1$ and $\mathcal{E}(\rho_0,\chi_{E_0}) <\infty$, we define the \textit{discrete flow}, or \textit{minimizing movements}, with time step size $h$ starting from $(\rho_0,\chi_{E_0})$ in the following way: $(\rho_0^h,\chi_{E^h_0}):=(\rho_0,\chi_{E_0})$ and, for every $n \in \mathbb{N} \setminus\{0\}$, we set
\begin{align}\label{minpb}
(\rho_n^h,\chi_{E_n^h})\in \arg\min_{(\rho, \chi_E)}	\Big\{  \frac{1}{2h} W^2_2(\rho, \rho^h_{n-1}) +  \mathcal{E}(\rho,\chi_E) : \rho \in K, \,\, E \text{ measurable}\Big\}.
\end{align}

The minimization problem~\eqref{minpb} is well defined since the existence of a minimizer follows from the direct method.
Indeed, from the convexity of $s \mapsto f(s)$, we have
\[
\int_\Omega f( \tilde \rho)\, \dx   \ge \int_\Omega f( \rho) \, \dx + \int_\Omega f'(\rho) (\tilde \rho-\rho)\, \dx ,
\]
for all functions $\rho$ and $\tilde \rho$.
Whence we deduce that $\mathcal{E}(\rho,\chi_E)$ is a lower semicontinuous functional with respect to $L^{m}$-weak convergence in $\rho$ and $L^1$-strong in $\chi_E$. The $L^{m}$-weak compactness of the variable $\rho$ follows from the Dunford-Pettis Theorem (see for instance~\cite[Box 8.2]{Santabook}) since  $\int_\Omega f(\rho)\le \mathcal{E}(\rho,\chi_E)$.
The $L^1$-strong compactness in the variable $\chi_E$ follows since $\int_{\Omega} 1 \, \mathrm{d}|\nabla \chi_E | \le \mathcal{E}(\rho,\chi_E)$.

We remark that the choice of the initial set $E_0$ is arbitrary, and we have the following uniform bound
\[
\int_\Omega 1 \, \mathrm{d} |\nabla \chi_{E^h_n}| \le \mathcal{E}(\rho_0,\chi_{E_0}).
\]

Let $(\rho^h_n,\chi_{E^h_n})$ be any minimizer of~\eqref{minpb}, by piecewise constant interpolation we define the \textit{approximate flow} $(\rho^h(t,\cdot),\chi_{E^h(t)})_{t \in [0,+\infty)}$ by
\begin{align} \label{eq:defdisc}
\rho^h(t,\cdot):=\rho_n^h(\cdot), \quad E^h(t):=E^h_n \quad \text{for every } t \in [nh,(n+1)h).
\end{align}
We will prove that, fixed $T>0$, the approximate flow $(\rho^h(t,\cdot),\chi_{E^h(t)})_{t \in [0,T)}$ converges as $h \to 0$ to a ``relaxed" distributional solution to the Verigin problem~\eqref{eq:smooth_system}. In particular, our main result reads as follows.

\begin{theorem}\label{thm:main}
Let $T>0$ and $\Omega \subset \R^d$ be an open set.    Let $E_0 \subset \Omega$ be a measurable set, and $\rho_0 \in L^1(\Omega)$ be initial conditions such that $\int_\Omega \rho_0=1$ and $\mathcal{E}(\rho_0,\chi_{E_0})< \infty$, and let $(\rho^h, E^h)$ be constructed as in~\eqref{minpb}-\eqref{eq:defdisc} for $h>0$. Then there exist a subsequence $h \rightarrow 0$, functions $\rho \in L^{m}(\Omega \times (0,T))$,  $\chi \in   L^1(\Omega \times (0,T))$ with $0\le \chi \le 1$, and a vector field $\mathbf{u} \in (L^1( \Omega \times (0,T); \rho \dx\dt ))^{ d}$ such that
    \begin{align}
    \rho^h &\rightarrow \rho \quad  \text{ strongly in } L^1(\Omega \times (0,T))  \text{ as } h \rightarrow 0, \\
    \rho^h &\rightharpoonup \rho \quad  \text{ weakly in } L^{m}(\Omega \times (0,T))  \text{ for } m>1  \text{ as } h \rightarrow 0, \\
    \chi_{E^h} &\stackrel{\ast}{\rightharpoonup} \chi  \quad  \text{ weakly-$*$ in } L^\infty(\Omega \times (0,T)) \text{ as } h \rightarrow 0,
    \end{align}
    and
    \begin{equation}\label{eq:transport}
        - \int_0^T \int_{\Omega} \rho (\partial_t \eta + \mathbf{u}  \cdot \nabla \eta ) \, \dx \dt = \int_{\Omega} \rho_0 \eta(\cdot, 0) \, \dx
    \end{equation}
    for all $\eta \in C^\infty_c(\overline\Omega \times [0, T))$.
      Moreover,  $(\rho, \chi) $ satisfies the energy dissipation rate 
    \begin{align}\label{eq:energydiss}
    \mathcal{E}(\rho(\tau), \chi(\tau)) + \int_0^\tau \int_\Omega \rho |\mathbf{u}|^2\, \dx \dt \leq \mathcal{E}(\rho_0, \chi_{E_0}),
    \end{align}
    for almost every $\tau\in (0,T)$.
   The triple $(\rho, \chi, \mathbf{u} )$ satisfies the distributional Young-Laplace law
    \begin{align}\label{eq:weaksoleq}
       & - \int_0^T \int_{\Omega} \rho \mathbf{u} \cdot \xi \,\mathrm{d} x \dt  \\
      \notag  &= \int_0^T \int_{\Omega} (\operatorname{div} \xi - \nu  \cdot \nabla \xi \nu  )
   \, \mathrm{d} |\nabla \chi| \dt 
      -\int_0^T \int_{\Omega} \Big(\rho^{m} - \lambda\chi\Big) \mathrm{div} \xi \,  \dx \dt,
    \end{align}
    for all $\xi \in C_c^\infty(\overline \Omega \times [0, T); \mathbb{R}^d) $ with $\xi \cdot \nu_{ \Omega} =0$ along $\partial \Omega \times (0,T)$, where $\nu $ denotes the Radon-Nikodym derivative $-\frac{ \mathrm{d} \nabla \chi}{\mathrm{d} |\nabla \chi|}$, and the optimal energy dissipation relation 
    \begin{align}\label{eq:optimalenrel}
&\mathcal{E}(\rho,\chi)+\frac{1}{2}\int_0^T
\int_{\Omega} \rho |\mathbf{u}|^2 \, \dx \dt -\frac{1}{2} \int_0^T \int_{\Omega} \rho |\xi|^2\, \dx \dt \\
&-\int_0^T\int_{\Omega} (\operatorname{div}\xi -\nu \cdot \nabla \xi \nu) \, \mathrm{d}|\nabla \chi| \dt + \int_0^T \int_\Omega \rho^m \operatorname{div} \xi \, \dx \dt -\int_0^T \int_{\Omega} \lambda\chi \operatorname{div} \xi  \, \dx \dt\notag \\
&
    \leq \mathcal{E}(\rho_0, \chi_{E_0}) \notag 
     \end{align}
     for all $\xi \in C_c^\infty(\overline \Omega \times [0, T); \mathbb{R}^d) $ with $\xi \cdot \nu_{ \Omega} =0$ along $\partial \Omega \times (0,T)$.
\end{theorem}

{Under the additional assumption that the functions $(\rho^h \wedge 1)^\frac{d}{d-1}$ are uniformly Muckenhoupt weights (cf. Definition \ref{def:Mucken}), we can improve the weak convergence of $(\chi_{E^h})_{h>0}$ to strong convergence, ensuring that the limit function $\chi$ is indeed the characteristic function of a set.}

\begin{theorem}
    \label{cor} 
    Let $\Omega = \R^d$ for $d\ge 3$.
	 Under the assumptions of Theorem \ref{thm:main} and the additional assumption that    $(\rho^h(\cdot, t) \wedge 1)^\frac{d}{d-1} \in A_1^*$ for every $h>0$ and $t \in [0,T)$, with $\sup_{h>0, t \in [0,T)]} c_{\rho^h(\cdot, t)}<\infty$, we have
	\begin{align} \label{eq:strongcompchi}
		\chi_{E^h} \rightarrow \chi  \quad  \text{ strongly in } L^1(\R^d \times (0,T); (\rho \wedge 1)^\frac{d}{d-1} \dx \dt ) \text{ as } h \rightarrow 0.
	\end{align}
    Moreover, there exists a one-parameter family of sets $(E(t))_{t \in [0,T)}$ of finite perimeter in $\Omega$ such that $\chi = \chi_E $ in $\{ \rho >0\}$.  
\end{theorem}

We observe that any sufficiently smooth triple $(\rho, \chi,\mathbf{u})$ satisfying \eqref{eq:transport} and~\eqref{eq:weaksoleq}, whose existence is established by the theorems above, is in fact a strong solution to the Verigin problem~\eqref{eq:smooth_system}.

\begin{proposition}
Let $(\rho, \chi,\mathbf{u})$ satisfy~\eqref{eq:transport} and~\eqref{eq:weaksoleq} and be such that $\rho \in C^\infty(\overline{\Omega} \times [0,T))$  and $\chi(\cdot,t)=\chi_{E(t)}$ for some family of smooth sets $(E(t))_{t \in [0,T)}$, then $(\rho, \chi,\mathbf{u})$ solves~\eqref{eq:smooth_system} with boundary conditions~\eqref{eq:boundcond1}-\eqref{eq:boundcond2}.
\end{proposition}

\begin{proof} 
    Integrating by parts~\eqref{eq:transport}, we obtain \begin{align}\label{eq:int_bp_transport}
       & \int_0^T \int_\Omega ( \partial_t \rho + [\![\rho]\!]V + \operatorname{div}(\rho \mathbf{u}))\eta \, \dx \dt \\ \notag
        &- \int_0^T \int_{\partial E} [\![\rho \mathbf{u}\cdot \nu_E]\!] \eta \, \dH \dt -\int_0^T\int_{\partial \Omega}  \rho  \mathbf{u}\cdot \nu_\Omega  \eta \, \dH \dt =0
    \end{align}
    for all $\eta \in C^\infty_c(\overline\Omega \times [0, T))$.
    Similarly, from~\eqref{eq:weaksoleq} it follows
    \begin{align*}
       - \int_0^T\int_\Omega \rho \mathbf{u}\cdot \xi \, \dx \dt
       = &\int_0^T \int_{\partial E} H  \nu_E \cdot \xi \, \dH \dt + \int_0^T \int_\Omega \nabla( \rho^m -\lambda\chi_E) \cdot \xi \, \dx \dt \\
       &-\int_0^T \int_{\partial E} [\![\rho^m-\lambda\chi]\!] \xi \cdot \nu_E \, \dH \dt
    \end{align*}
    for all $\xi \in C_c^\infty(\overline \Omega \times [0, T); \mathbb{R}^d)$ with $\xi \cdot \nu_{ \Omega} =0$ along $\partial \Omega \times (0,T)$. Since $\pi = \rho^m - \lambda \chi_E$, we can infer~\eqref{eq:smooth_system} with their natural boundary conditions~\eqref{eq:boundcond1} and~\eqref{eq:boundcond2}.
\end{proof}

\section{Time-discrete problem}\label{sec:discreteprob}
The following lemma provides the Euler--Lagrange equations satisfied by the discrete flow.

\begin{lemma}[Euler--Lagrange equations]\label{lemma:EL}
     Given~$h>0$, $E_0 \subset \Omega$  and $\rho_0 \in L^1(\Omega)$ such that $\int_\Omega \rho_0 \,\dx=1$ and $\mathcal{E}(\rho_0,\chi_{E_0}) <\infty$, and let $(\rho_n^h,\chi_{E_n^h})$  be constructed as in~\eqref{minpb}.
   Then, for any vector field $\xi\in C^1_c(\overline\Omega;\mathbb{R}^d)$ with $\xi \cdot \nu_\Omega =0$ on $\partial \Omega$, we have 
    \begin{align}\label{eq:eulerlagrdis}
       &- \int_{\Omega} \rho^h_n \nabla \phi^h_n \cdot \xi \,\mathrm{d} x   \\
      \notag  &= \int_{\partial^* E_n^h \cap \Omega} (\operatorname{div} \xi - \nu_{E^h_n} \cdot \nabla \xi \nu_{E^h_n} )
      \,\mathrm{d} \mathcal{H}^{n-1}-\int_{\Omega}  (\rho^h_n)^m \mathrm{div}\xi \,\mathrm{d} x
      +\int_{E^h_n} \lambda\mathrm{div}\xi \,\mathrm{d} x ,
    \end{align}
    where $ \nu_{E^h_n}$ denotes the outer normal vector field to $\partial^* E^h_n$.
\end{lemma}

\begin{proof}[Proof of Lemma~\ref{lemma:EL}]
Consider the flow map $\Phi_s:\Omega \to \Omega$ associated to $\xi$, i.e., $\partial_s\Phi_s=\xi \circ \Phi_s$ and $\Phi_0= \text{Id}$, where $\text{Id}$ denotes the identity map. We consider variations of $(\rho^h_n,\chi_{E^h_n})$ of the form $(\rho_s,\chi_{E_s})=((\Phi_s)_\#\rho^h_n, \chi_{E^h_n} \circ \Phi_s^{-1})$. 
Note carefully that this is indeed an admissible competitor: by standard ODE theory, since $\xi$ is tangential to $\partial \Omega$ by assumption, $\Phi_s$ leaves $\partial \Omega$ invariant, and hence also $\Omega$. So $\Phi_s\colon \Omega \to \Omega$ is a diffeomorphism for sufficiently small $s$.

Let $\phi^h_n$ be the Kantorovich potential such that $(\text{Id}-h \nabla \phi^h_n)_\# \rho^h_n = \rho^h_{n-1}$, using variations as above we obtain for the Wasserstein distance
\begin{align*}
    \frac{\mathrm{d}}{\mathrm{d} s}\Big |_{s=0} W_2^2(\rho_s,\rho^h_{n-1})
    &= \lim_{s \to 0} \frac{W_2^2(\rho_s,\rho^h_{n-1})-W^2_2(\rho^h_n,\rho^h_{n-1})}{s}.
\end{align*}
Choosing $\pi = (\Phi_s,\text{Id}+h \nabla \phi^h_n)_\# \rho^h_n$ as a competitor, we see that
\begin{align*}
W_2^2(\rho_s,\rho^h_{n-1})&=W_2^2((\Phi_s)_\#\rho^h_n,(\text{Id}-h \nabla \phi^h_n)_\# \rho^h_n)\\
&\le \int_\Omega |\Phi_s(x) -x+h \nabla \phi^h_n(x)|^2 \rho^h_n(x) \,\mathrm{d} x\\
&= \int_\Omega |s\xi(x) +o(s)+h \nabla \phi^h_n(x)|^2 \rho^h_n(x) \,\mathrm{d} x\\
&=W^2_2(\rho^h_n,\rho^h_{n-1})+2sh \int_{\Omega} \xi \cdot \nabla \phi^h_n \rho^h_n \,\mathrm{d} x
+o(s).
\end{align*}
Taking the limit as $s \to 0^+$ and $s \to 0^-$ we infer
\begin{align*}
    \frac{\mathrm{d}}{\mathrm{d} s}\Big |_{s=0} \frac{1}{2h}W_2^2(\rho_s,\rho^h_{n-1})=  \int_{\Omega} \xi \cdot \nabla \phi^h_n \rho^h_n \,\mathrm{d} x.
\end{align*}
The first variation of the perimeter term is
\[\frac{\mathrm{d}}{\mathrm{d} s}\Big |_{s=0} \int_{\Omega } 1 \, \mathrm{d} |\nabla \chi_{E_s}|= \int_{\partial^* E^h_n \cap \Omega}(\operatorname{div} \xi - \nu_{E^h_n} \cdot \nabla \xi \nu_{E^h_n} ) \,\mathrm{d} \mathcal{H}^{n-1},\]
cf.~\cite[Theorem 17.5]{Maggi}.
By the very definition of push-forward measure and the Monge-Amp\'ere equation $ (\rho_s \circ \Phi_s) J \Phi_s=\rho^h_n$ we get, for $m=1$,
\begin{align*}
    \frac{\mathrm{d}}{\mathrm{d} s}\Big |_{s=0} \int_{\Omega } \rho_s \log \rho_s \,\mathrm{d} x&= \frac{\mathrm{d}}{\mathrm{d} s}\Big |_{s=0} \int_{\Omega } \rho^h_n \log \rho_s \circ \Phi_s \,\mathrm{d} x\\
    &= -\frac{\mathrm{d}}{\mathrm{d} s}\Big |_{s=0} \int_{\Omega } \rho^h_n \log J\Phi_s \,\mathrm{d} x\\
    &=- \int_{\Omega } \rho^h_n \operatorname{div}\xi \,\mathrm{d} x,
\end{align*}
{whereas for $m>1$,}
{
\begin{align*}
    \frac{\mathrm{d}}{\mathrm{d} s}\Big |_{s=0} \int_{\Omega } \frac{1}{m-1}\rho_s^{m}  \mathrm{d} x&= \frac{\mathrm{d}}{\mathrm{d} s}\Big |_{s=0} \int_{\Omega } \frac{1}{m-1} \rho^h_n  (\rho_s \circ \Phi_s)^{m-1}\mathrm{d} x\\
    &= \frac{\mathrm{d}}{\mathrm{d} s}\Big |_{s=0} \int_{\Omega } \frac{1}{m-1}\rho^h_n \Big (\frac{\rho_n^h}{J\Phi_s} \Big )^{m-1}  \mathrm{d} x\\
    &= \frac{\mathrm{d}}{\mathrm{d} s}\Big |_{s=0} \int_{\Omega } \frac{1}{m-1}(\rho^h_n)^{m} (J\Phi_s^{-1})^{m-1}  \mathrm{d} x\\
    &=- \int_{\Omega } (\rho^h_n)^m \operatorname{div}\xi \,\mathrm{d} x.
\end{align*}
}
Using again the definition of push-forward and the definition of $E_s$, we notice
\begin{align*}
    -\frac{\mathrm{d}}{\mathrm{d} s}\Big |_{s=0} \int_{\Omega } \chi_{E_s} \rho_s \,\mathrm{d} x= -\frac{\mathrm{d}}{\mathrm{d} s}\Big |_{s=0} \int_{\Omega } \chi_{E_s}\circ \Phi_s \rho^h_n \,\mathrm{d} x = -\frac{\mathrm{d}}{\mathrm{d} s}\Big |_{s=0} \int_{\Omega } \chi_{E^h_n} \rho^h_n \,\mathrm{d} x=0.
\end{align*}
For the last term we get
\begin{align*}
    \frac{\mathrm{d}}{\mathrm{d} s}\Big |_{s=0} \lambda\int_{\Omega } \chi_{E_s}  \,\mathrm{d} x= \frac{\mathrm{d}}{\mathrm{d} s}\Big |_{s=0} \lambda\int_{\Omega } \chi_{E_s} \circ \Phi_s J\Phi_s  \,\mathrm{d} x = \lambda\int_{\Omega } \chi_{E^h_n} \operatorname{div} \xi  \,\mathrm{d} x.
\end{align*}
By combining the above calculations we obtain~\eqref{eq:eulerlagrdis}.
\end{proof}

The discrete flow constructed as in~\eqref{minpb}, in particular the sets $E_n^h$, satisfies the following important property, which will be fundamental to show the convergence of the perimeters of $E^h$ as $h \to 0$.

\begin{lemma}[Discrete almost minimizing property]
	Given~$h>0$, $E_0 \subset \Omega$  and $\rho_0 \in L^1(\Omega)$ such that $\mathcal{E}(\rho_0,\chi_{E_0}) <\infty$, let $(\rho_n^h,\chi_{E_n^h})$ be constructed as in~\eqref{minpb}. Then, it holds
\begin{align}\label{eq:almostEh}
	& \int_{\Omega} 1 \, \mathrm{d}|\nabla \chi_{E^h_n} |  
	- \int_{\Omega}  \chi_{E^h_n}  \rho^h_n \, \mathrm{d}x \\
	&\leq 
	\int_{\Omega} 1 \, \mathrm{d}|\nabla \chi_{E} | 
	- \int_{\Omega}  \chi_{E}  \rho^h_n  \, \mathrm{d}x 
	+ \lambda \int_{\Omega} | \chi_{E} - \chi_{E^h_n}|   \, \mathrm{d}x \notag .
\end{align}
\end{lemma}

\begin{proof}
	By testing the minimality of $(\rho_n^h,\chi_{E_n^h})$ with $(\rho^h_n, E)$, for an arbitrary $E \subset \Omega$, we obtain
	\begin{align*}
		\int_{\Omega} 1 \, \mathrm{d}|\nabla \chi_{E^h_n} |  
		\leq 
		\int_{\Omega} 1 \, \mathrm{d}|\nabla \chi_{E} |  
		+ \int_{\Omega} ( \chi_{E} - \chi_{E^h_n}) \Big( \lambda   - \rho^h_n\Big )  \, \mathrm{d}x ,
	\end{align*}
	whence we deduce \eqref{eq:almostEh}.
\end{proof}

\section{Compactness}\label{sec:compactness}
Using the energy dissipation of the scheme and standard Wasserstein estimates, we obtain the following weak compactness result for $\rho^h$.

\begin{lemma}
\label{lemma:compactness rho}
    Given~$h>0$, $E_0 \subset \Omega$  and $\rho_0 \in L^1(\Omega)$ such that $\int_\Omega \rho_0 \, \dx=1$ and $\mathcal{E}(\rho_0,\chi_{E_0}) <\infty$, and let $(\rho_n^h,\chi_{E_n^h})$ be constructed as in~\eqref{minpb}. Then
\begin{equation} \label{eq:discomp}
    W_2(\rho^h(t), \rho^h(s)) \leq \sqrt{2}\mathcal{E}(\rho_0, \chi_{E_0})^\frac12 (t-s)^\frac12,
\end{equation}
for $t>s \ge0$ with $t-s\ge h$.
Furthermore, for every $T>0$, we have
\begin{equation}\label{eq:weakcomprho}
\rho^h \rightharpoonup \rho \quad \text{ weakly in } L^{m}(\Omega \times (0,T))
\end{equation}
and
\begin{equation}\label{eq:contcomp}
    W_2(\rho(t), \rho(s)) \leq \sqrt{2}\mathcal{E}(\rho_0, \chi_{E_0})^\frac12 (t-s)^\frac12,
\end{equation}
for $t>s\ge 0$.
\end{lemma}
\begin{proof}
    Using $\rho^h_{n-1}$ as a competitor in~\eqref{minpb} yields
    \[
    \frac{1}{2h} W^2_2(\rho^h_n, \rho^h_{n-1}) + \mathcal{E}(\rho^h_n, \chi_{E^h_n}) \leq \mathcal{E}(\rho^h_{n-1}, \chi_{E^h_{n-1}}),
    \]
    whence, summing in $n$, we get 
    \begin{equation}\label{eq:sum}
    \frac{1}{2} \sum_{n=n_0+1}^{n_1}\frac{1}{h}W^2_2(\rho^h_n, \rho^h_{n-1}) + \mathcal{E}(\rho^h_{n_1}, \chi_{E^h_{n_1}}) \leq \mathcal{E}(\rho^h_{n_0}, \chi_{E^h_{n_0}}).
    \end{equation}
    For any pair of integers $n_1 > n_0 \geq 0$, we have
    \[
 \ W_2(\rho_{n_0}, \rho_{n_1}) \leq \sum_{n=n_0 +1}^{n_1}   W_2(\rho_{n}, \rho_{n-1})
    \]
due to the triangle inequality, and
    \[
  \Bigg(\frac{1}{n_1-n_0}\sum_{n=n_0 +1}^{n_1}   W_2(\rho_{n}, \rho_{n-1}) \Bigg)^2
 \leq   \frac{1}{n_1-n_0} \sum_{n=n_0 +1}^{n_1}   W_2^2(\rho_{n}, \rho_{n-1})
    \]
    due to the Jensen inequality.
    It follows from~\eqref{eq:sum} that
    \begin{align*}
        W_2(\rho^h_{n_1},\rho^h_{n_0}) &\leq \sqrt{(n_1 - n_0)h } \Bigg (\sum_{n=n_0+1}^{n_1}\frac{1}{h} W^2_2(\rho^h_n, \rho^h_{n-1})\Bigg)^{\frac{1}{2}} \\
        & \leq \sqrt{(n_1 - n_0)h } \sqrt{2(\mathcal{E}(\rho^h_{n_0},\chi_{E^h_{n_0}})-\mathcal{E}(\rho^h_{n_1},\chi_{E^h_{n_1}}))},
    \end{align*}
    which implies~\eqref{eq:discomp}.
    
   Since $\int_0^T \int_\Omega f(\rho^h) \, \dx \dt $   is uniformly bounded, the compactness~\eqref{eq:weakcomprho} follows directly for $m>1$, and 
    from~\eqref{eq:discomp} and Dunford-Pettis Theorem
     for $m=1$ (see for instance~\cite[Box 8.2]{Santabook}).
From the continuity of $W_2$ with respect to weak$*$-convergence (see~\cite[Theorem 5.10]{Santabook}), we deduce~\eqref{eq:contcomp}.
\end{proof}

We now turn to the weak compactness of the sets $E^h$.

\begin{lemma}
\label{lemma:compactness chi}
      Given~$h>0$, $E_0 \subset \Omega$  and $\rho_0 \in L^1(\Omega)$ such that $\int_\Omega \rho_0 \, \dx=1$ and $\mathcal{E}(\rho_0,\chi_{E_0}) <\infty$, and let $(\rho_n^h,\chi_{E_n^h})$ be constructed as in~\eqref{minpb}.
      Then,
      \begin{equation}\label{eq:weakcompchi}
\chi_{E^h} \stackrel{\ast}{\rightharpoonup} \chi \quad \text{ weakly-$*$ in } L^\infty(\Omega \times (0,\infty))
\end{equation}
and
\begin{equation}\label{eq:weakstarcomp Dchi}
\nabla \chi_{E^h}\stackrel{\ast}{\rightharpoonup} \nabla \chi \quad \text{ weakly-$*$ in } \mathcal{M}(\Omega\times(0,\infty);\mathbb{R}^d ).
\end{equation}
\end{lemma}
\begin{proof}[Proof of Lemma~\ref{lemma:compactness chi}]
    The convergence~\eqref{eq:weakcompchi} follows directly from the uniform boundedness of $\chi_{E^h}$ and the Banach-Alaoglu theorem.
    Since 
    \[
    \int_{\Omega \times (0,T)} 1 \,  \mathrm{d}|\nabla \chi_{E^h(t)}(x)|  \dt \le \int_0^T\mathcal{E}(\rho^h(t),\chi_{E^h}(t))\, \dt  \le T \mathcal{E}(\rho_0,\chi_{E_0})
    \]
    for every $h>0$, the Banach-Alaoglu theorem yields~\eqref{eq:weakstarcomp Dchi}.
\end{proof}

Furthermore, observe that, by Kantorovich duality we have
\begin{align}\label{eq:kanto}
\frac{1}{2h} W^2_2(\rho_n^h, \rho_{n-1}^h) = \frac{h}{2} \int_{\Omega} |\nabla \phi_n^h|^2 \rho_n^h \, \text{d}x,
\end{align}
where we recall that $\phi_n^h$ is the Kantorovich potential. 
Hence, choosing $n_0 = 0$ and $n_1=N$ in~\eqref{eq:sum},
we deduce 
\begin{equation*}
\begin{split}
    \frac{h}{2} \sum_{n=1}^{N} \int_{\Omega} |\nabla \phi_n^h|^2 \rho_n^h \, \text{d}x
    \le \mathcal{E}(\rho_{0},\chi_{E_{0}})-\mathcal{E}(\rho_{N}^h,\chi_{E_{N}^h}),
\end{split}
\end{equation*}
whence, letting $T=Nh$,
\begin{equation}\label{phibound}
    \frac{1}{2}\int_0^T \int_{\Omega} |\nabla \phi^h|^2 \rho^h \, \text{d}x \text{d}t \le \mathcal{E}(\rho_{0},\chi_{E_{0}}).
\end{equation}
This bound will be useful in the next sections.

An application of a variant of the the Aubin--Lions lemma allows us to strengthen the compactness properties of $\rho^h$.
\begin{proposition}\label{prop:strongrho}
    Given~$h>0$, $E_0 \subset \Omega$  and $\rho_0 \in L^1(\Omega)$ such that $\int_\Omega \rho_0 \, \dx=1$ and $\mathcal{E}(\rho_0,\chi_{E_0}) <\infty$, and let $(\rho_n^h,\chi_{E_n^h})$ be constructed as in~\eqref{minpb}. Then
    \begin{equation}\label{eq:strongcomprho}
    \rho^h \rightarrow \rho \quad \text{ strongly in } L^1(\Omega \times (0,T)). 
    \end{equation}
    Moreover, for $m>1$,
    \begin{align*}
    &\rho^h \rightarrow \rho \quad \text{ strongly in } L^p(\Omega \times (0,T)), \quad \text{ for } 1\le p<m.
    \end{align*}
\end{proposition}

To prove this proposition, we use the optimality condition which can be found in~\cite[Proposition 8.7]{Santabook}.
\begin{lemma}\label{phivarphi}
    For any $h>0$ and every $n=1,2,\ldots$, there exists a constant $c_n^h\in \mathbb{R}$ such that the optimal measure $\rho^h_n$ satisfies
        \begin{align*}
         f'(\rho^h_n) - \chi_{E_n^h} + \phi_n^h&= c_n^h  \quad \text{ a.e. in } \supp(\rho_n^h)
        \\
      f'(\rho^h_n) - \chi_{E_n^h} +\phi_n^h&\geq c_n^h \quad \text{ a.e. in } \Omega,
    \end{align*}
    where $f'(\rho^h_n) = \begin{cases}
        \log (\rho^h_n) +1, \; m=1,\\
        \frac{m}{m-1}(\rho^h_n)^{m-1}, \;  m>1,
    \end{cases}$ and $\phi^h_n$ is the (unique) Kantorovich potential from $\rho^h_n$ to $\rho_{n-1}^h$, i.e. $(\textnormal{Id}-h \nabla \phi^h_n)_\# \rho^h_n=\rho^h_{n-1}.$ 
\end{lemma}

\begin{proof}[Proof of Proposition~\ref{prop:strongrho}]

Since $\phi^h_n$ is Lipschitz continuous (see~\cite[Corollary 3.19]{AmbrosioBrueSemola}),
we have that  $f'(\rho^h_n) - \chi_{E_n^h}$ is Lipschitz. In particular, $\rho_n^h \in SBV(\Omega)$ and
{
\begin{align}
    \nabla^a f'(\rho^h_n) &= -\nabla \phi_n^h \quad \text{ a.e. in } \Omega, \label{eq:absgrad}\\
    [\![f'(\rho^h_n)]\!]&= f'(\rho^h_{n,+})-f'(\rho^h_{n,-})=-1 \quad \text{ on } \partial^* E_n^h \cap \Omega.
    \label{eq:jumplogrho}
\end{align}
In particular,  we infer $e \rho_{n,+}^h=\rho_{n,-}^h$ for $m=1$.}

Define an auxiliary function $g: (0, + \infty] \rightarrow (0, + \infty]$ as
{
\begin{align*}
    g(s) := \begin{cases}
        s^{m}, & s<1,\\
        \tilde g(s) , & s\ge 1,\\
        + \infty, & s=+\infty,
    \end{cases}
\end{align*}
where $\tilde g(s) := \begin{cases}
        \log (s) +1, \; m=1,\\
       s^{m-1}, \;  m>1,
    \end{cases}$ and
which is continuous, increasing, concave for $m=1$, and convex for $m>1$.
}

Consider $g\circ \rho^h_n$ and observe that, for every $n$,
\[
\int_\Omega |g (\rho^h_n)| \, \dx \leq C\mathcal{E}(\rho_{0},\chi_{E_{0}})+1.
\]
To estimate the total variation $|D(g(\rho_n^h))|(\Omega)$ we treat its absolutely continuous part and jump part separately. 
First, we compute
{
\begin{align*}
   \int_\Omega |\nabla^a (g \circ  \rho_n^h) | \, \dx &= 
    \int_{\{\rho <1\}} |\nabla^a   (\rho_n^h)^{m} |\, \dx
    + \int_{\{\rho \geq 1\}} |\nabla^a  \tilde g(\rho_n^h)|\, \dx \\
    &\leq  C \int_{\Omega} \sqrt{\rho_n^h}|\nabla^a   f'(\rho_n^h) |\, \dx,
\end{align*}}
whence, using~\eqref{eq:absgrad} and ~\eqref{phibound}, we obtain the uniform boundedness of
$  \int_0^T\int_\Omega |\nabla^a (g \circ  \rho^h) | \, \dx \dt$.
Then, to bound the jump part of $|D(g(\rho_n^h))|(\Omega)$, we consider three distinct cases
{
\begin{align*}
    |[\![g(\rho^h_n)]\!]| &= |(\rho_{n,+}^h)^{m}-(\rho_{n,-}^h)^{m}| < 1, \quad && \text{ for }\rho_{n,+}, \rho_{n,-} \le1,\\
    |[\![g(\rho^h_n)]\!]| &= |\tilde g(\rho_{n,+}^h)- \tilde g(\rho_{n,-}^h)| = 1, \quad && \text{ for }\rho_{n,+},  \rho_{n,-} \ge 1,\\
    |[\![g(\rho^h_n)]\!]| &= |\tilde g (\rho_{n,-}^h)-(\rho_{n,+}^h)^{m}| \leq c_m, \quad && \text{ for }\rho_{n,-} \ge 1 >\rho_{n,+}, 
\end{align*}
}
where we used~\eqref{eq:jumplogrho}, and $0<c_m \leq 3$.
It follows that
\begin{equation}
    \int_0^T \|(g \circ \rho^h)(\cdot,t)\|_{BV(\Omega)}\, \dt  \le C, 
\end{equation}
where $C>0$ is a uniform constant.

     We introduce the lower semicontinuous functional $\mathcal{F}: L^1(\Omega) \to [0,+\infty]$ defined as
        \[
    \mathcal{F}(u) := \int_\Omega f(u) \, \dx + \|g(u)\|_{BV(\Omega)}. 
    \]

Note that, for every $C>0$, the set $\mathcal{V}:=\{
v \in L^1(\Omega): \mathcal{F}(v) \le C\} 
$
is compact in $L^1(\Omega)$ with respect to the strong topology. Indeed, let $v_n$ be a sequence in $\mathcal{V}$, then, up to a subsequence, $g(v_n)  \rightarrow w$ strongly in $ L^1(\Omega)$ and almost everywhere in $\Omega$. Since $g$ is strictly monotone, thus invertible, 
$
v_n = g^{-1}(g(v_n)) \to g^{-1}(w)$ a.e. in $\Omega$.
 Moreover, by the dominated convergence theorem, $v_n  \rightarrow g^{-1}(w) $ in $L^1(\Omega)$.

Recalling also~\eqref{eq:discomp}, it follows that all the assumptions of Theorem~\ref{theo:RossiSava} are satisfied with $d=W_2$ (Wasserstein distance) and $\mathcal{U}=(\rho^h)_{h>0}$, in particular $\mathcal{U}$ is tight with respect to $\mathcal{F}$. Therefore, we apply Theorem~\ref{theo:RossiSava} and obtain that $\mathcal{U}$ is relatively compact in $\mathcal{M}(0,T; L^1(\Omega))$
. In particular, $\rho^h(\cdot,t) \rightarrow \rho(\cdot,t) $ strongly in $L^1(\Omega)$,  
for almost all $t \in (0,T)$. Finally, the strong convergence~\eqref{eq:strongcomprho} follows from an application of the dominated convergence theorem.
{Moreover, for $m>1$, the convergence holds in $L^{p}(\Omega)$ strongly for every $1\le p <m$ by interpolation.}
\end{proof}

\section{Proof of Theorem~\ref{thm:main}: Existence of distributional solutions}\label{sec:existence}

\subsection{Construction of the velocity field and transport equation~\eqref{eq:transport}}
Defining $\mathbf{u}^h := \nabla \phi^h$ and $\mathbf{v}^h := \sqrt{\rho^h} \nabla \phi^h $, we have (up to a subsequence) that 
\begin{equation} \label{eq:limitvh}
\sqrt{\rho^h} \mathbf{u}^h = \mathbf{v}^h \rightharpoonup \mathbf{v} \quad \text{weakly in } L^2(\Omega \times (0,T)).
\end{equation}
Furthermore, one can show that the measure $\sqrt{\rho}\mathbf{v}  \mathcal{L}^d \llcorner \Omega \otimes \mathcal{L}^1 \llcorner (0,T)$  is absolutely continuous with respect to $ \rho  \mathcal{L}^d \llcorner \Omega \otimes \mathcal{L}^1 \llcorner (0,T)$, in particular there exists a vector field $\mathbf{u} \in L^1( \Omega \times (0,T); \rho \dx\dt )$ such that $\sqrt{\rho}\mathbf{v}  =  \rho\mathbf{u}  $ in $L^1(\Omega \times (0,T))$ and
\begin{equation}\label{eq:identificationv}
    \mathbf{v}  = \sqrt{\rho}\mathbf{u}  \quad \text{on }\supp \rho.
\end{equation}
Indeed, consider two open sets $\mathcal{U} \subset \Omega $ and $I \subset (0, T)$, and $ \xi \in  C^\infty_c(\overline \Omega\times (0,T)) $ such that $\|\xi\|_\infty\leq 1$ and $\supp \xi \subset \mathcal{U} \times I$, then the convergences~\eqref{eq:strongcomprho} and~\eqref{eq:limitvh}, Young's inequality, and the estimate~\eqref{phibound} yield
\begin{align*}
   &\int_0^T \int_{\Omega}  \mathbf{v}  \cdot \xi \sqrt{\rho}\, \dx \dt \\&= \lim_{h \rightarrow 0} 
    \int_0^T \int_{\Omega} \mathbf{v}^h  \cdot \xi \sqrt{\rho^h}\, \dx \dt \\
   & \leq  \bigg(\limsup_{h \rightarrow 0}  \int_0^T \int_{\Omega} | \nabla \phi^h|^2  \rho^h\, \dx \dt 
   \bigg)^\frac12
   \bigg(\limsup_{h \rightarrow 0}  \int_I \int_{\mathcal{U}} \rho^h| \xi|^2 \, \dx \dt 
   \bigg)^\frac12
   \\
   & \leq  \large(2 \mathcal{E}(\rho_{0},\chi_{E_{0}})
   \large)^\frac12
   \bigg(  \int_I \int_{\mathcal{U}} \rho \, \dx \dt 
   \bigg)^\frac12.
\end{align*}

Recalling that $(\text{Id}-h \nabla \phi^h_n)_\# \rho^h_n = \rho^h_{n-1}$, for $\eta \in C_c^\infty(\overline\Omega \times [0, +\infty))$ it holds
\begin{align*}
\frac1h \int_\Omega (\rho^h_n - \rho^h_{n-1}) \eta \, \dx = \frac1h \int_\Omega \rho^h_n(x) (\eta(x)- \eta(x-h \nabla \phi^h_n(x)) \, \dx .
\end{align*}
By Taylor's expansion in the right hand side, 
\begin{align}\label{eq:equationrho}
\frac1h \int_\Omega (\rho^h_n - \rho^h_{n-1}) \eta \, \dx =  \int_\Omega \rho^h_n(x) \nabla \eta \cdot   \nabla \phi^h_n  \, \dx  + \mathcal{R},
\end{align}
where the error $\mathcal{R}$ can be bounded by
\[
\frac{h}2  \|\nabla^2 \eta \|_\infty \int_\Omega |\nabla \phi^h_n |^2 \rho^h_n\, \dx = \frac{1}{2h} \|\nabla^2 \eta \|_\infty  W^2_2(\rho_n^h, \rho_{n-1}^h).
\]
Recalling~\eqref{eq:sum}, we obtain
\[
h \sum_{n=1}^N \frac{1}{2h} \|\nabla^2 \eta \|_\infty  W^2_2(\rho_n^h, \rho_{n-1}^h) \leq h \|\nabla^2 \eta \|_\infty \mathcal{E}(\rho_{0},\chi_{E_{0}}),
\]
in particular the error term $\mathcal{R}$ vanishes as $h \rightarrow 0$.
Hence, integrating in time~\eqref{eq:equationrho}, using~\eqref{eq:strongcomprho},~\eqref{eq:limitvh} and~\eqref{eq:identificationv}, we obtain~\eqref{eq:transport} in the limit $h \rightarrow 0$.

\subsection{Perimeter convergence}
We pass to the limit $h \to 0$ in the discrete almost minimizing property for $E^h$, given by~\eqref{eq:almostEh}, obtaining
\begin{align}\label{eq:almostminprop}
		\int_{0}^{T}\int_{\Omega} 1 \, \mathrm{d} |\nabla \chi| \dt -  \int_{0}^{T}\int_{\Omega} \chi(\rho -\lambda) \, \dx \dt \leq 	\int_{0}^{T}\int_{\Omega} 1 \, \mathrm{d} |\nabla \tilde \chi| \dt -  \int_{0}^{T}\int_{\Omega} \tilde \chi(\rho -\lambda)\, \dx \dt 
	\end{align}
for any $\tilde \chi : \Omega \to \{0,1\}$ measurable, where we recall $\chi$ is the weak$^*$-limit of $\chi_{E^h}$.

	We know that $(\chi_{E^h}, \rho^h)$ satisfies~\eqref{eq:almostEh}, namely
	\begin{align} \label{eq:halmostminprop}
	&\int_{0}^{T}	\int_{\Omega} 1 \, \mathrm{d} |\nabla \chi_{E^h}|\dt  - \int_{0}^{T}\int_{\Omega} \chi_{E^h}(\rho^h -\lambda) \, \dx \dt
    \notag \\
    &\leq \int_{0}^{T}	\int_{\Omega} 1 \, \mathrm{d} |\nabla \tilde \chi|\dt  - \int_{0}^{T}\int_{\Omega} \tilde \chi(\rho^h -\lambda)\, \dx \dt,
	\end{align}
	for any $\tilde \chi : \Omega  \to \{0,1\}$ measurable.
	Note that, since $\tilde \chi,\, \chi_{E^h}  \in L^\infty(\Omega)$ and $\rho^h \rightarrow \rho$ strongly in $L^1(\Omega \times (0,T))$ (cf.~\eqref{eq:strongcomprho}), we have
	\[
	\lim_{h \to 0}  \int_{0}^{T}\int_{\Omega} \tilde \chi\rho^h \, \dx \dt =   \int_{0}^{T} \int_{\Omega} \tilde \chi\rho \, \dx \dt 
	\]
    and
    \begin{align}\label{eq:claim}
	\lim_{h \to 0}  \int_{0}^{T} \int_{\Omega}  \chi_{E^h}\rho^h \, \dx \dt =    \int_{\Omega} \chi\rho \, \dx \dt .
\end{align}
	Moreover, using~\eqref{eq:weakstarcomp Dchi}, by the lower semi-continuity of the total variation, we have
	\begin{equation}\label{eq:semicontchi}
		\int_{0}^{T}\int_{\Omega} 1 \, \mathrm{d} |\nabla \chi| \dt \leq \liminf_{h \to 0} \int_{0}^{T}\int_{\Omega} 1 \, \mathrm{d} |\nabla \chi_{E^h}| \dt .
	\end{equation}
Combining these together we infer~\eqref{eq:almostminprop}.

We observe that~\eqref{eq:halmostminprop} holds also for $f : \Omega \rightarrow [0,1]$, instead of $\tilde \chi : \Omega \rightarrow \{0,1\}$.
Indeed, let $f : \Omega \rightarrow [0,1]$ and let $\chi_s := 1_{\{f >s\}}$. Using $\tilde \chi = \chi_s$ in~\eqref{eq:halmostminprop}, integrating over $s \in (0,1)$, the coarea and layer-cake formulas yield
\begin{align*} 
&	\int_{0}^{T}	\int_{\Omega} 1 \, \mathrm{d} |\nabla \chi_{E^h}|\dt  - \int_{0}^{T}\int_{\Omega} \chi_{E^h}(\rho^h -\lambda) \, \dx \dt\\
	&\leq \int_0^1   \Big(\int_{0}^{T}	\int_{\Omega} 1 \, \mathrm{d} |\nabla  \chi_s|\dt  - \int_{0}^{T}\int_{\Omega} \chi_s(\rho^h -\lambda)\, \dx \dt \Big) \, \mathrm{d} s 
\\
& \leq \int_{0}^{T}	\int_{\Omega} |\nabla f| \, \mathrm{d} x \dt  - \int_{0}^{T}\int_{\Omega} f(\rho^h -\lambda)\, \dx \dt .
\end{align*}

Similarly, using $\tilde \chi = \chi $ (i.e., the weak limit of $\chi_{E^h}$ as $h \rightarrow 0$ given by~\eqref{eq:weakcompchi}) in~\eqref{eq:halmostminprop}, then taking the limit $h \rightarrow 0$, we obtain
\begin{align*} 
	&	\limsup_{h \to 0} \int_{0}^{T}	\int_{\Omega} 1 \, \mathrm{d} |\nabla \chi_{E^h}|\dt 
	\\
	&\leq \int_{0}^{T}	\int_{\Omega} 1 \, \mathrm{d} |\nabla  \chi|\dt + \lim_{h \to 0} \int_{0}^{T}\int_{\Omega}( \chi_{E^h}- \chi)\rho^h  \, \dx \dt ,
\end{align*}
where the last term is zero due to~\eqref{eq:weakcompchi} 
and~\eqref{eq:strongcomprho}.

This combined with the lower semi-continuity~\eqref{eq:semicontchi} gives
\begin{equation}\label{eq:convperim}
 \lim_{h \to 0} \int_{0}^{T}\int_{\Omega} 1 \, \mathrm{d} |\nabla \chi_{E^h}| \dt = \int_{0}^{T}\int_{\Omega} 1 \, \mathrm{d} |\nabla \chi| \dt .
\end{equation}
In particular, the above limit implies
\begin{equation}\label{eq:weakstartv}
    |\nabla \chi_{E^h}|\stackrel{\ast}{\rightharpoonup} |\nabla \chi| \quad \text{ weakly-$*$ in } \mathcal{M}(\Omega\times(0,+\infty) ).
\end{equation}

\subsection{Energy dissipation rate~\eqref{eq:energydiss}}
\label{sec:endissrate}

Similarly as in~\cite{Chambolle-Laux}, we can show 
that the De Giorgi variational interpolation
\begin{align*}
\big(\tilde\rho^h, \tilde{E}^h\big)((n-1)h +t) \in  \arg \min_{(\rho,E)} \Big\{ \frac{1}{2t} W^2_2(\rho, \rho^h_{n-1}) + \mathcal{E}(\rho, \chi_{E}) \Big\}
\end{align*}
satisfies 
\begin{equation}\label{dgine}
\begin{split}
    &\frac{h}{2} \left (\frac{W_2(\rho_n^h, \rho_{n-1}^h)}{h} \right)^2 + \frac{1}{2} \int_{0}^{h} \left (\frac{W_2(\tilde \rho^h((n-1)h+t),\rho_{n-1}^h)}{t} \right )^2 \text{d}t \\
    &\le \mathcal{E}(\rho_{n-1}^h,\chi_{E_{n-1}^h})-\mathcal{E}(\rho_{n}^h,\chi_{E_{n}^h}).
\end{split}
\end{equation}
Recall that $\tilde \rho^h ((n-1)h + \cdot ) = \rho^h_{n}$ and 
$\tilde E^h ((n-1)h + \cdot ) = E_{n}^h$
for any $t \in (0,h]$.

Indeed, w.l.o.g. by assuming $n=1$ and omitting the index $h$, we define 
\begin{equation*}
f(t):= \frac{1}{2t} W^2_2(\tilde\rho(t), \rho_{0}) + \mathcal{E}(\tilde\rho(t), \chi_{\tilde{E}(t)})
\end{equation*}
to be the minimal value in the problem above. Using the minimality of $(\tilde\rho, \tilde{E})$, we compute for $s < t$
\begin{align*}
f(t) - f(s) \leq \frac{s-t}{2st} W^2_2(\tilde\rho(s), \rho_{0}), 
\end{align*}
which yields
\begin{align*}
    \frac{f(t) - f(s)}{t-s} \leq - \frac{1}{2st} W^2_2(\tilde\rho(s), \rho_{0}) \rightarrow - \frac{1}{2t^2} W^2_2(\tilde\rho(t), \rho_{0}) \quad \text{as } s \rightarrow t. 
\end{align*}
By interchanging the roles of $s$ and $t$, i.e. using $s > t$, we can obtain the analogous reverse inequality. Hence, we have that $f$ is locally Lipschitz on $(0,h]$ with 
\begin{align}
    \frac{\mathrm{d}}{\dt} f(t) = - \frac{1}{2t^2} W^2_2(\tilde\rho(t), \rho_{0})
\end{align}
for almost every $t \in (0,h)$. For $\varepsilon>0$, we can integrate this equality from $t=\varepsilon$ to $t=h$, obtaining
\begin{align}
& \mathcal{E}(\tilde\rho(h), \chi_{\tilde{E}(h)}) + \frac{1}{2h} W^2_2(\tilde\rho(h), \rho_{0})  + \int_\varepsilon^h \frac{1}{2t^2} W^2_2(\tilde\rho(t), \rho_{0}) \, \dt  \\
 &=  \mathcal{E}(\tilde\rho(\varepsilon), \chi_{\tilde{E}(\varepsilon)}) + \frac{1}{2\varepsilon} W^2_2(\tilde\rho(\varepsilon), \rho_{0}).  \notag
\end{align}
Note that the right hand side is bounded by $\mathcal{E}(\rho_0, \chi_{{E}_0})$ due to the minimizing movements scheme. Hence,  for $n=1$ the equation~\eqref{dgine} follows by noticing that $\tilde\rho(0)= \rho_{0}$, $\tilde E(0)=E_0$, $\tilde \rho(h)=\rho^h_1$, and $\tilde E(h)=E^h_1$.

Summing in~\eqref{dgine} over $n$ from $n_0 +1$ to $n_1$, we obtain
\begin{equation}
\begin{split}
    &\frac{h}{2} \sum_{n=n_0+1}^{n_1} \left (\frac{W_2(\rho_n^h, \rho_{n-1}^h)}{h} \right)^2 + \frac{1}{2} \sum_{n=n_0+1}^{n_1}  \int_{0}^{h} \left (\frac{W_2(\tilde \rho^h((n-1)h+t),\rho_{n-1}^h)}{t} \right )^2 \text{d}t \\
    &\le \mathcal{E}(\rho_{n_0}^h,\chi_{E_{n_0}^h})
    - \mathcal{E}(\rho_{n_1}^h,\chi_{E_{n_1}^h}) \le \mathcal{E}(\rho_{0},\chi_{E_{0}})
    - \mathcal{E}(\rho_{n_1}^h,\chi_{E_{n_1}^h}).
    \label{eq:discdiss}
\end{split}
\end{equation}
First, recalling~\eqref{eq:kanto}-\eqref{phibound} and, in particular, the convergence~\eqref{eq:limitvh}, we can pass to the limit $h \to 0$ in the first term on the left hand side of~\eqref{eq:discdiss}.
Then, note that, by a change of variable,
\[
 \frac{1}{2}\sum_{n=n_0+1}^{n_1} \int_{0}^{h} \left (\frac{W_2(\tilde \rho^h((n-1)h+t),\rho_{n-1}^h)}{t} \right )^2 \, \text{d}t =  \frac{1}{2}\int_{n_0 h}^{n_1 h} \left (\frac{W_2(\tilde \rho^h(t),\rho^h(t))}{t-[t/h]h} \right )^2 \,  \text{d}t ,
\]
and
\[
\frac{W^2_2(\tilde \rho^h(t),\rho^h(t))}{2(t-[t/h]h)}  = \frac{t-[t/h]h}{2} \int_{\Omega} |\nabla \tilde \phi^h(t)|^2 \tilde \rho^h \, \dx,
\]
where $\tilde \phi^h$ is the optimal Kantorovich potential in $W_2(\tilde \rho^h(t),\rho^h(t))$.
For $n_0 = 0$, $n_1=N$, and $\tau=Nh \in (0,T)$, it follows that 
\[
\frac{1}{2}\sum_{n=+1}^{N} \int_{0}^{h} \left (\frac{W_2(\tilde \rho^h((n-1)h+t),\rho_{n-1}^h)}{t} \right )^2 \text{d}t = \frac{1}{2} \int_{0}^{\tau} \int_{\Omega} |\nabla \tilde \phi^h(t)|^2 \tilde \rho^h \,  \dt\dx.
\]
Recalling that $\tilde \rho^h((n-1)h+\cdot)=\rho_n^h$, we deduce that also $\mathbf{\tilde v}^h := \sqrt{\tilde \rho^h} \nabla \tilde \phi^h $ weakly converges in $L^2(\Omega \times (0,T))$ to the function $\mathbf{v}$ given by~\eqref{eq:limitvh}, hence we can also pass to the limit in the second term in the left hand side of~\eqref{eq:discdiss}. Moreover, observe that
\begin{align*}
\int_0^{\tau} \int_\Omega |\mathbf{v}|^2\, \dx \dt 
\geq \int_{\supp \rho} |\mathbf{v}|^2\, \dx \dt = \int_{\supp \rho} \rho |\mathbf{u}|^2\, \dx \dt 
= \int_0^{\tau} \int_\Omega \rho |\mathbf{u}|^2\, \dx \dt,
\end{align*}
where we used~\eqref{eq:identificationv}.

To conclude we show the lower semicontinuity of the energies $\mathcal{E}(\rho^h(\cdot,t),\chi_{E^h(t)})$ for almost every $t \in (0,T)$.
Since $\chi_{E^h}$ and $|\nabla \chi_{E^h}|$ only converge weakly-* in space-time, we first establish the lower semicontinuity against any non-negative test function $\xi \in C_c((0,T))$. The desired result then follows from the arbitrariness of $\xi$.

By Helly's selection theorem, up to a subsequence, $\mathcal{E}(\rho^h(\cdot,t),\chi_{E^h(t)}) \to e(t)$ as $h \to 0$ for almost every $t \in (0,T).$ 
In particular,  for every $\xi \in C_c( (0,T))$, by the dominated convergence theorem, it holds
\begin{equation}\label{eq:convHelly}
    \lim_{h \to 0} \int_0^T \xi(t) \mathcal{E}(\rho^h(\cdot,t),\chi_{E^h(t)}) \, \dt= \int_0^T \xi(t) e(t)\, \dt.
\end{equation}

Let $\xi \in C_c( (0,T))$ be such that $\xi \ge 0$, we claim that
\begin{equation}\label{eq:lsc_test}
    \liminf_{h \to 0} \int_0^T \xi(t) \mathcal{E}(\rho^h(\cdot,t),\chi_{E^h(t)}) \, \dt \ge \int_0^T \xi(t) \mathcal{E}(\rho(\cdot,t),\chi(\cdot,t)) \, \dt.
\end{equation}
Using~\eqref{eq:strongcomprho} combined with~\eqref{eq:strongcompchi} 
we deduce
\[
\lim_{h \to 0}  \int_0^T \xi(t) \int_\Omega \chi_{E^h(t)}(x)\rho^h(x,t)\, \dx \dt=\int_0^T\xi(t)\int_\Omega \chi(x,t) \rho(x,t)\, \dx \dt.
\]
Applying the weak-$*$ convergence~\eqref{eq:weakcompchi} with test function $\xi$,  we obtain
\[
\lim_{h \to 0} \int_0^T \xi(t) \int_\Omega \chi_{E^h(t)}(x)\, \dx \dt=\int_0^T \xi(t)\int_\Omega \chi(x,t)\, \dx \dt.
\]
Similarly, from~\eqref{eq:weakstartv} using as test function $t \mapsto \xi(t)$ multiplied by a suitable cut-off in the $x$ variable, we infer
\[
\lim_{h \to 0} \int_0^T \xi(t) \int_\Omega |\nabla \chi_{E^h(t)}(x)|\, \dx \dt \ge \int_0^T \xi(t)\int_\Omega |\nabla \chi(x,t)|\, \dx \dt.
\]
Finally, combining everything together and by the lower semicontinuity of  $\int_\Omega f(\rho^h) \, \dx$, the claim~\eqref{eq:lsc_test} follows. 

By the arbitrariness of $\xi$ and~\eqref{eq:convHelly}, we deduce that $e(t) \ge \mathcal{E}(\rho(\cdot,t),\chi(\cdot,t))$ for almost every $t \in (0,T)$, which implies the lower semicontinuity of $\mathcal{E}(\rho^h(\cdot,t),\chi_{E^h(t)})$.

\subsection{Passage to the limit in the Euler-Lagrange equation}\label{sec:limitEL}
From the convergences~\eqref{eq:strongcomprho} and~\eqref{eq:limitvh} it follows that
\begin{equation*}
    \rho^h \nabla \phi^h \to \sqrt{\rho}\textbf{v} \quad \text{ weakly in } L^{1}(\Omega \times (0,T)) ,
\end{equation*}
as $ h \to 0$.
Hence, the left hand side of~\eqref{eq:eulerlagrdis} passes to the desired limit as $h \rightarrow 0$.
Moreover, we can pass to the limit in the second and third terms in the right hand side of~\eqref{eq:eulerlagrdis} thanks to~\eqref{eq:weakcomprho} and~\eqref{eq:weakcompchi}, respectively.
Finally, the convergence of the first term of the right hand side of~\eqref{eq:eulerlagrdis} as $h \rightarrow 0$ follows from~\eqref{eq:convperim} and an application of Reshetnyak's continuity theorem~\cite[Thereom 2.39]{AFP}, and it reads as
\[
\int_0^T\int_{\Omega} (\operatorname{div} \xi - \nu_{E^h_n} \cdot \nabla \xi \nu_{E^h_n} )
 \, 
 \mathrm{d} |\nabla \chi_{E_n^h}| \dt  \rightarrow 
 \int_0^T\int_{\Omega} (\operatorname{div} \xi - \nu \cdot \nabla \xi \nu )
 \, 
 \mathrm{d} |\nabla \chi| \dt
\]
as $h \rightarrow 0$,
where $\nu$ denotes the Radon-Nikodym derivative $-\frac{\mathrm{d} \nabla \chi}{\mathrm{d} |\nabla \chi|}$.

\subsection{Optimal energy dissipation relation~\eqref{eq:optimalenrel}}

The local slope of $\mathcal{E}$, defined by
\[
|\partial \mathcal{E} (\bar \rho, \bar \chi)| :=\limsup_{\substack{\rho \to \bar \rho \\   \chi  \to \bar \chi   }} \frac{(\mathcal{E}(\bar \rho,\bar \chi)-\mathcal{E}(\rho, \chi))_+}{W_2(\bar \rho,\rho)},
\]
satisfies (cf.~\cite[Lemma 3.1.3]{AGS})
\[
|\partial \mathcal{E}|(\tilde \rho^h((n-1)h+t), \tilde \chi_{E^h((n-1)h+t)})\le \frac{W_2(\tilde \rho^h((n-1)h+t),\rho^h_{n-1})}{t}.
\]
Plugging this into~\eqref{eq:discdiss} yields the inequality
\begin{equation}\label{eq:disenslope}
\begin{split}
    &\frac{h}{2} \sum_{n=n_0+1}^{n_1} \left (\frac{W_2(\rho^h_n,\rho^h_{n-1})}{h} \right )^2 +\frac{1}{2}\int_{n_0 h}^{n_1 h} |\partial \mathcal{E}|^2 (\tilde \rho^h(t),\chi_{\tilde E^h(t)}) \, \dt \\
    &\le 
    \mathcal{E}(\rho_{n_0}^h,\chi_{E_{n_0}^h})
    - \mathcal{E}(\rho_{n_1}^h,\chi_{E_{n_1}^h}).
\end{split}
\end{equation}

Given any smooth divergence free vector field $\xi \in C^\infty(\Omega; \mathbb{R}^d)$, we obtain the one-parameter family $\rho_s$ via the continuity equation $\partial_s \rho_s + \xi \cdot \nabla \rho_s =0$.
Similarly as in~\cite[Section 4, Step 4]{Chambolle-Laux}, we can deduce 
\[
\frac12 |\partial \mathcal{E}|^2(\rho,\chi_E) \geq \lim_{s \to 0} \frac1s (\mathcal{E}(\rho,\chi_E) - \mathcal{E}(\rho_s,\chi_{E_s}))_+ - \frac1{2s^2}W^2_2(\rho, \rho_s). 
\]
Since $\xi$ is divergence free, by 
means of the change of variable $s' \mapsto s's$, we obtain the continuity equation $\partial_{s'}\rho_{s'}+\text{div}(s\xi \rho_{s'})=0$, in particular $\rho_{s'} |_{s'=1 }= \rho_s$.
Hence, using the Benamou-Brenier formula for the Wasserstein distance, we obtain
\[
W^2_2(\rho,\rho_s) \le \int_0^1 \int_{\Omega} \rho_{s'} |s\xi|^2 \,\dx \mathrm{d}s' = s^2 \int_{\Omega} \rho |\xi|^2 \,\dx+o(s^2).
\]
On the other hand, we have{
\begin{align*}
    (\mathcal{E}(\rho,\chi_E)-\mathcal{E}(\rho_s,\chi_{E_s}))_+ 
    &\ge \mathcal{E}(\rho,\chi_E)-\mathcal{E}(\rho_s,\chi_{E_s})\\
    &=-s \frac{\mathrm{d}}{\mathrm{d}s}\Big|_{s=0} \mathcal{E}(\rho_s,\chi_{E_s}) + o(s)\\
    &=-s\int_{\Omega} (\operatorname{div} \xi - \nu_{E} \cdot \nabla \xi \nu_{E} )\,
  \mathrm{d} |\nabla \chi_E|
  \\
  & \quad +s\int_{\Omega} \rho^{m} \mathrm{div}\xi \,\mathrm{d} x
  -s\int_{E} \lambda\mathrm{div}\xi \,\mathrm{d} x +o(s),
\end{align*}
}
where the last equality follows from~\eqref{eq:eulerlagrdis}.
The inequalities above yield
\begin{align*}
    \frac{1}{2}|\partial \mathcal{E}|^2(\rho,\chi_E)&\ge -\int_{\Omega} (\operatorname{div} \xi - \nu_{E} \cdot \nabla \xi \nu_{E} )\,
  \mathrm{d} |\nabla \chi_E|\\
  &\quad  +\int_{\Omega} \rho^m \mathrm{div}\xi \,\mathrm{d} x
  -\int_{E} \lambda\mathrm{div}\xi \,\mathrm{d} x -\frac{1}{2}\int_\Omega \rho |\xi|^2 \, \dx.
\end{align*}
Combining the above inequality for $\rho=\rho^h$ and $\chi_E=\chi_{E^h}$ with~\eqref{eq:disenslope}, then passing to the limit $h \to 0$, we obtain~\eqref{eq:optimalenrel} by exploiting~\eqref{eq:kanto}-\eqref{eq:limitvh}, the arguments of Section~\ref{sec:limitEL}, and those at the end of Section~\ref{sec:endissrate}.

\section{Proof of Theorem~\ref{cor}:  Strong compactness for $\chi_{E^h}$}\label{sec:strongconv}
In this section we show a strong compactness property for $\chi_{E^h}$.
To this aim, we work under the technical assumptions $\Omega=\R^d$ and that $\rho^h$ are uniformly in the class of Muckenhoupt weights $A_1^\ast$ (see Appendix \ref{sec:Mucken}), which are needed to apply the weighted Sobolev isoperimetric inequality.

The proof of Theorem~\ref{cor} relies on the following crucial stability estimate that allows us to transfer time-compactness from $\rho^h$ to $\chi_{E^h}$ and states in particular that the nonlinear splitting~$\rho \mapsto (\phi,\chi)$ into a smooth ($H^1(\rho)$) function and a pure jump function according to
\[
	f'(\rho) - \chi  + \phi = \text{const.} \quad \text{$\rho$-a.e.\ in } \R^d
\]
is unique and stable. 
Note that here, $\phi$ will play the role of the Kantorovich potential when transporting one time step to the previous one. 

\begin{lemma}\label{lemma:compchit}
	Let $(\rho,\chi,\phi)$ and $(\tilde \rho, \tilde \chi, \tilde \phi)$ where $\chi,\tilde \chi \colon \R^d \to \{0,1\}$ and $\rho, \tilde \rho$ are probability measures on $\Omega$ be such that
	{\begin{align}
		\label{eq:phirhochi}  f'(\rho) - \chi + \phi&= c  \quad \rho\text{-a.e.}
		\\
		\label{eq:phirhochiTILDE} f'(\tilde \rho) - \tilde \chi +\tilde \phi&= \tilde c  \quad \tilde\rho\text{-a.e.}
	\end{align}}
	for two constants $c,\tilde c\in \mathbb{R}$.
	Assume that $(\rho\wedge 1)^{\frac{d}{d-1}} , \, (\tilde \rho\wedge 1)^{\frac{d}{d-1}} \in A^*_1$. Assume also that there exists $C'>0$ such that $\|\rho\|_{L^{m}(\R^d)},\,\|\tilde \rho\|_{L^{m}(\R^d)},$ $\int_{\R^d} 1\, \mathrm{d}|\nabla \chi| $, $\int_{\R^d} 1\,\mathrm{d}|\nabla \tilde\chi|  \le C'$.
	Then there exists a constant $C=C(d,c_{\rho},c_{\tilde \rho}, C')>0$ such that
		\begin{align}
			&  \int_{\R^d} |\chi - \tilde\chi| (\rho\wedge \tilde \rho \wedge 1)^{\frac{d}{d-1}} \, \dx\notag\\
			&\leq C \|\rho-\tilde \rho\|_{L^p(\R^d)}  +  C\|\rho-\tilde \rho\|_{L^p(\R^d)}^{\frac{d}{2(d-1)}}\bigg(  \int_{\R^d} |\nabla \phi|^2 \rho \,\dx +  \int_{\R^d} |\nabla \tilde \phi|^2 \tilde \rho \,\dx\bigg)^{\frac{d}{2(d-1)}}, \label{eq:chi stability}
		\end{align}
        where $p=1$ for $1 \leq m \leq 2$ and $p= \frac{m}2$ for $m>2$.
\end{lemma}

\begin{proof}[Proof of Lemma~\ref{lemma:compchit}]
	As~$\eqref{eq:chi stability}$ only depends on~$\phi$ and~$\tilde \phi $ through their gradients, upon redefining them to~$\phi+c$ and~$\tilde \phi + \tilde c$, we may, without loss of generality, assume that~$c=\tilde c=0$.
	
    Since $\chi,\,\tilde \chi \in \{0,1\}$, using~\eqref{eq:phirhochi}--\eqref{eq:phirhochiTILDE}, we have
		\begin{equation}\label{stime}
        \begin{split}
			|\chi-\tilde \chi| 
			&=(\chi-\tilde \chi)( f'(\rho)- f'(\tilde \rho) + \phi - \tilde \phi ) \\
			&\le  |f'(\rho)- f'(\tilde \rho)| + (\chi-\tilde \chi)( \phi-\tilde \phi ).
        \end{split}
		\end{equation}

	To control the first right-hand side term in \eqref{stime} we use the local Lipschitz estimate
\begin{equation}\label{eq:estimate_mge1}
			|f'(\rho)- f'( \tilde \rho)| 
			\leq c_m \rho^{m-2} \vee  \tilde \rho^{m-2}|\rho -\tilde \rho|  ,
		\end{equation} 
        for some constant $c_m>0$ depending only on $m$.
Taking the integral of the first term of~\eqref{stime} multiplyed by $(\rho \wedge \tilde \rho \wedge 1)^\beta$, $\beta \in \{1, \frac{d}{d-1} \}$, and using \eqref{eq:estimate_mge1}, we obtain
		\begin{equation*}
        \begin{split}
			&\int_{\R^d} |f'(\rho)- f'( \tilde \rho)| (\rho \wedge \tilde \rho \wedge 1 )^\beta \, \dx \le c_m
 \int_{\R^d}  \rho^{m-2} \vee  \tilde \rho^{m-2}|\rho -\tilde \rho| (\rho \wedge \tilde \rho \wedge 1)^\beta\, \dx .
        \end{split}
		\end{equation*}
Note that for $1 \leq m\leq 2$ we have  $(\rho^{m-2} \vee  \tilde \rho^{m-2}) (\rho \wedge \tilde \rho \wedge 1)^\beta \le 1$, hence
        \begin{equation}\label{eq:estimate_mge2}
            \int_{\R^d} |f'(\rho)- f'( \tilde \rho)| (\rho \wedge \tilde \rho \wedge 1)^\beta \, \dx \le c_m    \|\rho-\tilde \rho\|_{L^{1}(\R^d)}.
        \end{equation}
On the other hand, for  $m > 2$, an application of H\"older's inequality gives
        \begin{equation}\label{eq:estimate_mle2}
        \begin{split}
            &\int_{\R^d} |f'(\rho)- f'( \tilde \rho)| ( \rho \wedge \tilde \rho \wedge 1)^\beta\, \dx\\
            &\le \int_{\R^d} |f'(\rho)- f'( \tilde \rho)|  \, \dx\\
            &\le C \max \big \{ \|\rho\|_{L^m(\R^d)}^{(m-2)/m},\|\tilde\rho\|_{L^m(\R^d)}^{(m-2)/m} \big\} \|\rho-\tilde \rho\|_{L^{m/2}(\R^d)}.
        \end{split}
        \end{equation}

	To estimate the last term on the right-hand side of~\eqref{stime}, we define
	\[g(u):= \begin{cases}
		|u|^{\frac{2d}{d-1}} & |u|\le 1,\\
		(2|u|-1)^{\frac{d}{d-1}} & |u| \ge 1.
	\end{cases}\]
	Let $g^*:\mathbb{R} \to \mathbb{R}$ be the convex conjugate of $g$, i.e., $g^*(v):=\sup \{u v- g(u): u \in \mathbb{R} \}$.
	It is easy to see that $g^*(0)=0$ and that $g^*$ is an even function with $\alpha_d:=g^*(1)=g^*(-1)<\infty$. In fact, one can compute
		\[
		\alpha_d= 
		\frac{d+1}{2d} \Big( \frac{d-1}{2d}\Big)^{\frac{d-1}{d+1}}
		\in (0,1)\quad \text{for any } d\geq 2.
		\]
	Therefore, using Fenchel's inequality and $\chi-\tilde \chi \in \{-1,0,1\}$
	\begin{align*}
		(\chi-\tilde\chi)( \phi-\tilde \phi) 
		&\leq   g^*(\chi-\tilde\chi)+  g(\phi-\tilde \phi)
		= \alpha_d  |\chi-\tilde\chi| +  g(\phi-\tilde \phi).
	\end{align*}
Since $\alpha_d<1$, we can absorb the first term into the left-hand side. To handle the second one, we construct a non-negative even $C^2$ function
	$\Phi: \mathbb{R} \to \mathbb{R}$ such that for some $C_\Phi>0$, we have $|u|^2 \wedge |u| \leq C_\Phi \Phi(u)$, $\Phi'(\{-1,0,1\})=\{0\}$,   $|\Phi'(u)|  \leq C_\Phi (1+\sqrt{|u|})$, and $|\Phi''(u)|\le C_{\Phi}$ (see, e.g., Figure \ref{Figura}). In particular $\Phi'$ is uniformly $\frac{1}{2}$-H\"older continuous in $\R$.

\begin{figure} \centering
\begin{tikzpicture}
  \begin{axis}[
    width=10cm, height=5cm,
    xlabel={$u$},  ylabel={$\Phi(u)$},
    xmin=-2.9, xmax=2.9,
    ymin=-0.05, ymax=0.4,
    axis x line=middle, 
    axis y line=none, 
    samples=800,
    tick align=outside,
    xtick=\empty,
    ytick=\empty,
    xlabel style={at={(axis description cs:0.95,0.1)},anchor=north east},
  ]
    \addplot[thick,domain=-2.3:2.3] 
      { (abs(x) <= 1) 
         ? ( (abs(x))^2 - (4/3)*(abs(x))^3 + 0.5*(abs(x))^4 ) 
         : ( (abs(x)) - 17/6 + (abs(x)+1)*exp(-(abs(x)-1)) ) 
      };
    \draw[dashed,gray!60] (axis cs:-1,0.0) -- (axis cs:-1,0.7);
    \draw[dashed,gray!60] (axis cs: 0,0.0) -- (axis cs: 0,0.7);
    \draw[dashed,gray!60] (axis cs:  1,0.0) -- (axis cs:  1,0.7);
    \node at (axis cs:-1,-0.06) [anchor=south,font=\footnotesize] {$-1$};
    \node at (axis cs: 0,-0.06) [anchor=south,font=\footnotesize] {$0$};
    \node at (axis cs: 1,-0.06) [anchor=south,font=\footnotesize] {$+1$};
    \node at (axis cs: 2.4,0.2) [anchor=south,font=\footnotesize] {\normalsize $\Phi(u)$};
    \draw[->,>=stealth] (axis cs:-2.8,-0.02) -- (axis cs:-2.8,0.4);
  \end{axis}
\end{tikzpicture}
\caption{Example of $\Phi$} \label{Figura}
\end{figure}
    
	Using the properties of the function $\Phi$, Theorem \ref{thm:wsii} with $w=\rho^{\frac{d}{d-1}}$, the Cauchy-Schwarz inequality, and then Young's inequality, we obtain 
	\begin{align*}
		\bigg ( &\int_{\R^d} g(\phi-\tilde\phi) (\rho \wedge \tilde \rho \wedge 1)^\frac{d}{d-1} \, \dx  \bigg)^{\frac{d-1}{d}} \\
		&\leq
		C_\Phi\bigg (\int_{\R^d} \Phi(\phi-\tilde\phi)^{\frac{d}{d-1}}  (\rho \wedge \tilde \rho \wedge 1)^\frac{d}{d-1}\, \dx \bigg)^{\frac{d-1}{d}} \\
		&\leq C_\Phi C_{Iso}
		\int_{\R^d} \big|\nabla \big(\Phi\circ(\phi-\tilde\phi) \big) \big| \rho \wedge \tilde \rho \wedge 1 \,\dx \\
		& \leq C
		\bigg (\int_{\R^d} |
		\Phi'(\phi-\tilde \phi)|^2  \, \rho \wedge \tilde \rho \wedge 1 \,\dx \bigg)^{\frac{1}{2}} 
		\bigg (\int_{\R^d} |
		\nabla (\phi-\tilde \phi)|^2 \, \rho \wedge \tilde \rho \wedge 1 \,\dx \bigg)^{\frac{1}{2}}\\
        & \leq C
		\bigg (\int_{\R^d} |
		\Phi'(\phi-\tilde \phi)|^2  \rho \wedge \tilde \rho \wedge 1\, \dx \bigg)^{\frac{1}{2}} 
		\bigg (\int_{\R^d} |
		\nabla \phi|^2  \rho\, \dx +
        \int_{\R^d} |
		\nabla \tilde \phi|^2  \tilde \rho\, \dx \bigg)^{\frac{1}{2}}.
	\end{align*}
    Moreover, since  $\chi,\tilde \chi \in \{0,1\}$ implies $\Phi'(\chi-\tilde \chi)=0$ and $\Phi'$ is $\frac{1}{2}$-H\"older, we have 
    	\begin{align*}
		|\Phi'(\phi-\tilde \phi)|^2
		&= 
		|\Phi'(\phi-\tilde \phi)-\Phi'(\chi-\tilde \chi )|^2\\
        &\leq C_\Phi \big|\phi - \tilde \phi -\chi +\tilde \chi\big|\\
        &\le C_\Phi  \big|f'(\rho)-f'(\tilde \rho)\big|,
	\end{align*}
	where we used again the equations~\eqref{eq:phirhochi}--\eqref{eq:phirhochiTILDE}.
	{Hence, we conclude using~\eqref{eq:estimate_mge2} and \eqref{eq:estimate_mle2}.}
\end{proof}

\begin{proof}[Proof of Theorem \ref{cor}.]
	Let $\tau, t \in (0,T)$.
	Thanks to Lemma~\ref{phivarphi},
	we can apply Lemma~\ref{lemma:compchit} to $( \rho,  \chi, \phi) = (\rho^h(t),\chi_{E^h}(t), \phi^h(t))$ and $(\tilde \rho ,\tilde \chi, \tilde \phi) = (\rho^h( t +\tau),\chi_{E^h(t+\tau)}, \phi^h(t+\tau))$. By using the estimates
	$(\rho(t) \wedge 1)^\frac{d}{d-1} \lesssim |\rho(t) \wedge 1 -\rho^h(t) \wedge 1|^\frac{d}{d-1} + (\rho^h(t) \wedge 1)^\frac{d}{d-1}$
	and
	$(\rho(t) \wedge 1)^\frac{d}{d-1} \lesssim |\rho(t) \wedge 1-\rho^h(t) \wedge 1|^\frac{d}{d-1} + |\rho^h(t) \wedge 1-\rho^h(t+\tau) \wedge 1|^\frac{d}{d-1} + (\rho^h(t+\tau) \wedge 1)^\frac{d}{d-1}$, and by integrating in $t$, we obtain
    \begin{align*} 
		&  \int_0^{T-\tau}\int_{\R^d} |\chi_{E^h(t+\tau)}(x)-\chi_{E^h(t)}(x)| (\rho(x,t) \wedge 1)^\frac{d}{d-1} \, \dx \dt \notag\\
		&\leq
		C\int_0^{T-\tau}
        \int_{\R^d}   |\rho(x,t) \wedge 1-\rho^h(x,t) \wedge 1|^\frac{d}{d-1} \, \dx \dt\\
        &\quad +
		C \int_0^{T-\tau}\int_{\R^d}   |\rho^h(x,t+\tau) \wedge 1-\rho^h(x,t) \wedge 1|^\frac{d}{d-1} \, \dx \dt   \\
         &\quad + C 
		 \int_0^{T-\tau} \Big(\int_{\R^d}   |\rho^h(x,t+\tau) -\rho^h(x,t) |^p \, \dx \Big)^\frac{1}{p} \dt   \\
        &  \quad +  C \int_0^{T-\tau} \bigg(\int_{\R^d}   |\rho^h(x,t+\tau) -\rho^h(x,t) |^p \, \dx \bigg)^{\frac{d}{2(d-1)p}} \\
        & \hspace{1cm} \cdot \bigg(  \int_{\R^d} |\nabla \phi^h(x,t)|^2 \rho^h(x,t) \,\dx +  \int_{\R^d} |\nabla \phi^h(x, t+\tau)|^2  \rho^h(x, t+\tau)\,\dx\bigg)^{\frac{d}{2(d-1)}} \hspace{-0.05cm}\dt .
    \notag 
	\end{align*} 
    For $d >2$, let $q'= \frac{2(d-1)}{d-2}$ be the conjugate exponent of $q=\frac{2(d-1)}{d} >1 $, then we can apply H\"older's inequality and the estimate \eqref{phibound} to bound the last term above, yielding
\begin{align*}
        & 
        \bigg[ \int_0^{T-\tau} \bigg(\int_{\R^d}   |\rho^h(x,t+\tau) -\rho^h(x,t) |^p \, \dx \bigg)^{\frac{1}{qp}} \\
        & \hspace{1cm}\cdot  \bigg(  \int_{\R^d} |\nabla \phi^h(x,t)|^2 \rho^h(x,t) \,\dx +  \int_{\R^d} |\nabla \phi^h(x, t+\tau)|^2  \rho^h(x, t+\tau)\,\dx\bigg)^{\frac{1}{q}} \dt \bigg]
        \\
        &\leq 
        \bigg[ \int_0^{T-\tau} \Big (\int_{\R^d}   |\rho^h(x,t+\tau) -\rho^h(x,t) |^p \, \dx \Big )^{\frac{d}{(d-2)p}} \dt \bigg]^\frac{1}{q'} \\
        & \quad \cdot \bigg[ \int_0^{T-\tau} \hspace{-0.09cm}\bigg(  \int_{\R^d} |\nabla \phi^h(x,t)|^2 \rho^h(x,t) \,\dx +  \int_{\R^d} |\nabla \phi^h(x, t+\tau)|^2  \rho^h(x, t+\tau)\,\dx\bigg) \dt \bigg]^{\frac{1}{q}}
        \\
        &\leq C  \bigg(\sup_t   \Big(\int_{\R^d}   |\rho^h(x,t+\tau) -\rho^h(x,t) |^p \, \dx \Big)^{\frac{2}{(d-2)p}} \bigg)^\frac{1}{q'}  \\
        & \quad  \cdot
        \bigg[ \int_0^{T-\tau} \Big (\int_{\R^d}   |\rho^h(x,t+\tau) -\rho^h(x,t) |^p \, \dx \Big )^{\frac{1}{p}} \dt \bigg]^\frac{1}{q'}    (2\mathcal{E}(\rho_{0},\chi_{E_{0}}) )^{\frac{1}{q}}\\
        &\leq C (T-\tau)^{\frac{p-1}{q'p}} \bigg(\sup_t   \Big(\int_{\R^d}   |\rho^h(x,t+\tau) -\rho^h(x,t) |^p \, \dx \Big)^{\frac{1}{(d-1)p}} \bigg)\\
        & \quad  \cdot
        \bigg[ \int_0^{T-\tau} \Big (\int_{\R^d}   |\rho^h(x,t+\tau) -\rho^h(x,t) |^p \, \dx \Big ) \dt \bigg]^\frac{1}{q'p}    (2\mathcal{E}(\rho_{0},\chi_{E_{0}}) )^{\frac{1}{q}},
    \end{align*}
    where in the last line we used the converse Jensen inequality.

	For $z\in \Omega$ with $|z|$ small, 
	we obtain the  space compactness
	\begin{align*}
		&\int_0^T\int_{\Omega \cap (\Omega-z)} |\chi_{E^h(t)}(x+z)-\chi_{E^h(t)}(x)| (\rho(x,t) \wedge 1)^\frac{d}{d-1} \, \dx \dt 
		\\&\leq C
		\int_0^T\int_{\Omega }  |\rho(x,t) - \rho^h(x,t)|^\frac{d}{d-1}\, \dx \dt 
		\\& \quad + C\int_0^T\int_{\Omega \cap (\Omega-z)} |\chi_{E^h(t)}(x+z)-\chi_{E^h(t)}(x)| (\rho^h(x,t) \wedge 1)^\frac{d}{d-1} \, \dx \dt 
		\\&\leq \int_0^T\int_{\Omega }  |\rho(x,t) - \rho^h(x,t)|^\frac{d}{d-1}\, \dx \dt   + T \sup_{t} \bigg( |z|   \int_\Omega 1 \, \mathrm{d} |\nabla \chi_{E^h(t)}|
		\bigg)
		\\
		&   \leq \int_0^T\int_{\Omega }  |\rho(x,t) - \rho^h(x,t)|\, \dx \dt 
		+T|z|  \mathcal{E}(\rho_0,\chi_{E_0}).
	\end{align*}
	This gives the strong convergence~\eqref{eq:strongcompchi} thanks to the strong convergence~\eqref{eq:strongcomprho} of $\rho^h$ and an application of the Frech\'et-Kolmogorov theorem.

Moreover, there exists a measureable set $E=\bigcup_{t\in(0,T)} E(t) \times\{t\} \subset \R^d\times(0,T)$ such that $E(t)$ is a set of finite perimeter for a.e.\ $t\in(0,T)$ and $\chi=\chi_E$.
Indeed, by the coarea formula and Fubini's theorem
\[
    \int_0^T \int_{\R^d} 1 \, \mathrm{d}|\nabla \chi(\cdot,t)|\, \dt =  \int_0^1 \Big (\int_0^T P (\{\chi(\cdot,t)>s\})\,\dt \Big ) \mathrm{d} s.
\]
Hence, there exists $s\in (0,1)$ such that $E:=\{\chi>s\}\subset \R^d\times(0,T)$ satisfies $\int_0^T P(E(t))\, \dt \leq \int_0^T \int_{\R^d} 1 \, \mathrm{d} |\nabla \chi(\cdot,t)|\, \dt <\infty$. Since $\chi\in \{0,1\}$ a.e. in $\{\rho>0\}$  due to~\eqref{eq:strongcompchi}, we have $\chi_E = \chi$ in $\{\rho>0\}$.
\end{proof}

\appendix

\section{Compactness in measure}

In this section, we state a variant of the Aubin--Lions
Lemma, which can be found in~\cite[Theorem 2]{Rossi-Savare}.

Let $B$ be a separable Banach space. We say that $\mathcal{F}:B \to [0,+\infty]$ is a 
\emph{coercive integrand} if
\[
\{v \in B : \mathcal{F}(v) \le C\} \text{ is compact for every } C \ge 0.
\] 

We denote by $\mathcal{M}(0,T;B)$ the space of Bochner-measurable $B$-valued functions with the topology induced by the convergence in measure, i.e. $(u_n)_{n \in \mathbb{N}} \subset \mathcal{M}(0,T;B)$ converges in measure to $u \in \mathcal{M}(0,T;B)$ as $n \to \infty$ if
\[
\lim_{n \to \infty} |\{t \in (0,T) :\|u_n(t)-u(t)\|_{B} \ge \sigma \}| = 0 \,\, \forall \sigma >0.
\]

\begin{theorem}\label{theo:RossiSava}
    Let $\mathcal{U}$ be a family of measurable $B$-valued functions; if there exist a lower semicontinuous coercive integrand $\mathcal{F}: B \to [0,+\infty]$ and a distance $d:B \times B \to [0,+\infty]$, such that
    \[
    \mathcal{U} \text{ is tight w.r.t } \mathcal{F}, \text{ i.e. } \quad S:=\sup_{u \in \mathcal{U}} \int_0^T \mathcal{F}(u(t)) \, \dt < \infty,
    \]
    and
    \[
    \lim_{\tau \to 0} \sup_{u \in \mathcal{U}} \int_0^{T-\tau} d(u(t+\tau),u(t))\, \dt =0,
    \]
    then $\mathcal{U}$ is relatively compact in $\mathcal{M}(0,T;B)$.
\end{theorem}

\section{Muckenhoupt weights}\label{sec:Mucken}
We introduce the notion of \emph{Muckenhoupt weight}, for a complete treatise on the argument we refer the interested reader to~\cite[Section 15]{NLPT} for weighted Sobolev spaces and~\cite{Baldi} for the weighted $BV$ space.

\begin{definition}[Muckenhoupt weights]\label{def:Mucken}
Let $\Omega \subset \R^d$ be an open set (not necessary bounded). Let $p\in (1,+\infty)$ and let $w \in L^1_{loc}(\Omega)$ be a nonnegative function such that $0<w<+\infty$ almost everywhere. We say that $w$ belongs to the \emph{Muckenhoupt class} $A_p$ if there exists a constant $c_{p,w}>0$ such that
\begin{equation*}
    \sup \Big ( \fint_B w \, \dx \Big ) \Big ( \fint_B w^{1/(1-p)} \, \dx \Big )^{p-1} =c_{p,w} <\infty,
\end{equation*}
where the supremum is taken over all balls $B \subset \Omega$. 
We say that $w \in A_1$ if for every ball $B\subset \Omega$ it holds
\begin{equation*}
    \fint_B w \, \dx \le c_w \essinf_B w,
\end{equation*}
for some constant $c_w>0$. Finally, if $w\in A_1$ is also lower semi-continuous we say that it belongs to $A_1^*$.
\end{definition}

We state a version of the Sobolev isoperimetric inequality that holds in the weighted $BV$ space.
\begin{theorem}{\cite[Theorem 1.3]{BGKM}}\label{thm:wsii}
   Let $w \in A_1^*$. Then there exists $C>0$ only depending on $c_w$ and on the dimension $d$ such that, for every $f \in BV(\R^d;w)$, it holds
   \[
   \|f\|_{L^{1^*}(\R^d;w)} \le C \|Df\|_{w^{1/1^*}}(\R^d).
   \]
\end{theorem}

\section*{Acknowledgements}
The authors thank Aymeric Baradat and Inwon Kim for helpful discussions and feedback on the problem.  Part of this contribution was completed while A.K.\ and A.M.\ were visiting T.L.\ at the University of Regensburg, whose hospitality and financial support through the Research Training Group RTG 2339 IntComSin is acknowledged. A.K.\ and A.M.\ thank the Centre International de Rencontres Math\'ematiques (CIRM) for hosting and supporting part of this work by the Research in Residence program.
A.K.\ research has been supported by the Austrian Science Fund (FWF) through grants 10.55776/F65, 10.55776/P35359, 10.55776/Y1292. 
A.M.\ research is sponsored by the Alexander von Humboldt Foundation, which is acknowledged. A.M.\ is funded by the Deutsche Forschungsgemeinschaft (DFG,
German Research Foundation) - CRC 1720 - 539309657.

\bibliographystyle{abbrv}
\bibliography{biblio}

\end{document}